\documentclass{cmslatex}
\usepackage[paperwidth=7in, paperheight=10in, margin=.875in]{geometry}
 \usepackage[backref,colorlinks,linkcolor=red,anchorcolor=green,citecolor=blue]{hyperref}
\usepackage{amsfonts,amssymb}
\usepackage{amsmath}
\usepackage{graphicx}
\usepackage{cite}
\usepackage{enumerate}
\usepackage{subfigure}
\usepackage{epsfig}
\usepackage{caption}
\usepackage{epstopdf}
\renewcommand{\theequation}{\arabic{section}.\arabic{equation}}

   \allowdisplaybreaks
\begin{document}
 \title{Interaction of the elementary waves for shallow water equations with discontinuous topography}

          %For each author, make a block with the following macros:

          \author{Qinglong Zhang\thanks{School of Mathematics and Statistics, Ningbo University,
Ningbo, 315211, P.R.China, (zhangqinglong@nbu.edu.cn).}
          \and Wancheng Sheng\thanks{Department of Mathematics, Shanghai University,
Shanghai, 200444, P.R.China, (mathwcsheng@t.shu.edu.cn).}
          \and Yuxi Zheng\thanks{Department of Mathematics, The Pennsylvania State University,
University Park, 16801, U.S.A, (yuz2@psu.edu.cn).}}

         \pagestyle{myheadings} \markboth{INTERACTION OF THE ELEMENTARY WAVES OF SHALLOW WATER FLOW}{ZHANG,SHENG AND ZHENG} \maketitle

% ----------------------------------------------------------------
\begin{abstract}
The Riemann problem of one dimensional shallow water equations with discontinuous topography has been constructed recently. The elementary waves include shock waves, rarefaction waves, and the stationary wave. The stationary wave appears when the water depth changes, especially when there exists a bottom step. In this paper, we are mainly concerned with the interaction between a stationary wave with either a shock wave or a rarefaction wave. By using the characteristic analysis methods, the evolution of waves is described during the interaction process. The solution in large time scale is also presented in each case. The results may contribute to research on more complicated wave interaction problems.

\

\noindent%
{\sc Keywords.}~ Shallow water equations; source term; interaction of elementary waves;
Riemann problem.

\vskip 6pt
\noindent%
{\sc 2010 AMS Subject Classification.} Primary: 35L65, 35L80, 35R35, 35L60; Secondary: 35L50.
\end{abstract}
%%%%%%%%%%%%%%%%%%%%%%%%%%%%%%%%%%%%%%%%%%%%%%%%%%%%%%

\maketitle
% ----------------------------------------------------------------

\section{Introduction}\label{sec1}
The one-dimensional (1D) shallow water equations are given by
\begin{equation}\label{1.1}
\left\{
\begin{array}{l}
h_t+(h u)_{x}=0,\\[5pt]
\displaystyle (h u)_t+(h(u^2+\frac{gh}2))_x=-gha_x,\\[5pt]
a_t=0, \\
\end{array}
\right.
\end{equation}
where $h$ denotes the height of the water from the bottom to the surface, $u$ is the velocity of the fluid, $g$ the gravity constant, and $a$ the height of the river bottom from a given level. See Fig. 1.1. System \eqref{1.1} expresses the nonlinear shallow water fluids with a step like bottom. While the system is equivalent to the isentropic gas flow with the adiabatic exponent $\gamma=2$, the derivative of \eqref{1.1} is totally different from gas dynamics since the height $y=h(x,t)$ is a free boundary which imposes free boundary conditions to the system. See \cite{Stoker} for more details about the derivation of \eqref{1.1}.

The shallow water model is widely used for free-surface flows arising in shores, rivers and other physical phenomena \cite{Stoker,Toro}. It is well known that system \eqref{1.1} is not conservative because the existence of source term $-gha_x$, i.e., it belongs to the class of resonance systems. The theory of nonconservative system is introduced in \cite{DalLeflochMurat4} and the follow-up works \cite{Goatin5,Isaacson}. For the other similar models, we refer the reader to 1D fluids in a variable cross-section duct \cite{Andrianov1,LeflochThanh12,Thanh17,Zhang} and multiphase flow models \cite{Andrianov2,SaurelAbgrall15}, to name just a few. 

\begin{figure}[h]
\centering
\includegraphics[width=0.56\textwidth]{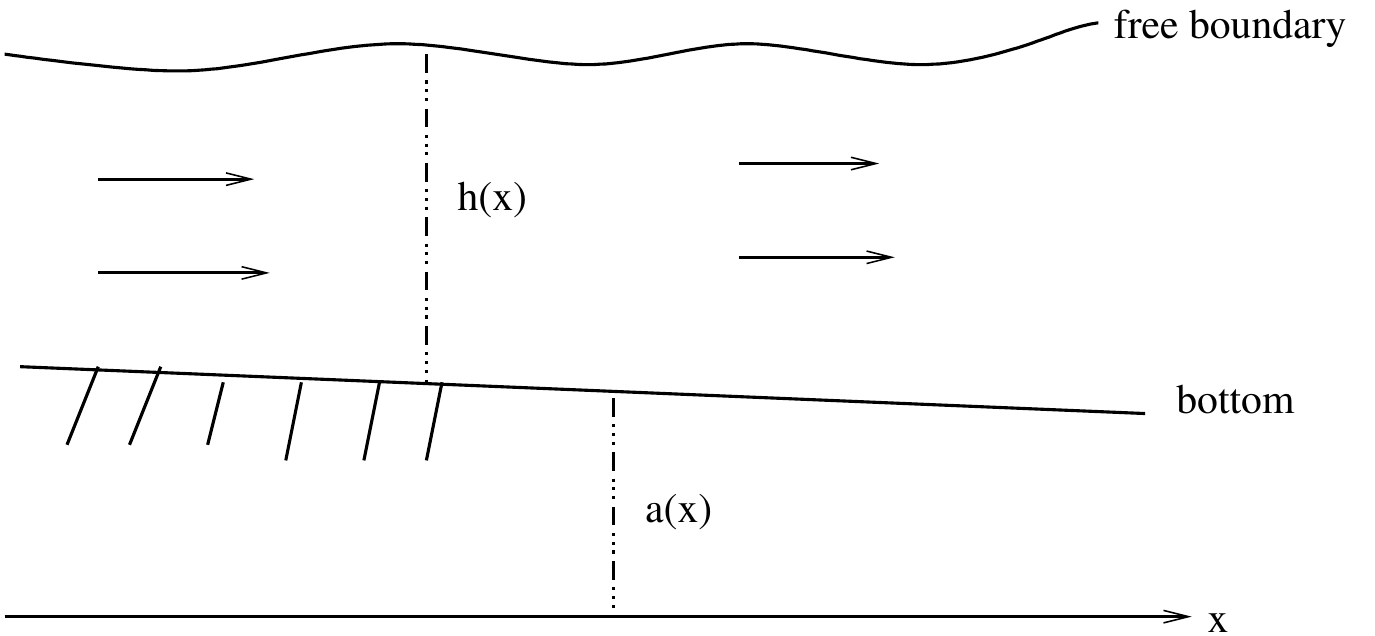}
\caption*{Fig. 1.1.The model of shallow water equation \eqref{1.1}.}
\end{figure}

The Riemann problem for \eqref{1.1} is given by
\begin{equation}\label{1.2}
(u,h,a)(x,0)=\left\{
\begin{array}{ll}
(u_-,h_-,a_0),\quad x<0,\\[5pt]
(u_+,h_+,a_1),\quad x>0.
\end{array}
\right.
\end{equation}
In  2003, LeFloch and Thanh (\cite{LeflochThanh12}) solved the Riemann problem of flows with a variable cross-section area. By dividing the phase  plane into several areas using certain hypersurfaces, the non-strict system can be viewed as strictly hyperbolic in each area. For the Riemann problem of the shallow water equations, they applied a similar method and constructed the solutions in the $(u,h)$ phase plane \cite{LeflochThanh13}. Alcrudo and Benkhaldoun \cite{Alcrudo} solved the Riemann problem of \eqref{1.1} with \eqref{1.2} from another approach, similar result can be found in \cite{Bernetti}. Besides, numerical approaches for system \eqref{1.1} are also available. Greenberg and Leroux  \cite{Greenberg} studied the discretization of source terms in nonlinear hyperbolic balance laws. Jin and Wen \cite{Jin1,Jin2} proposed interface-type well-balanced methods to capture steady state solutions for hyperbolic systems with geometric source terms. LeFloch and Thanh \cite{LeflochThanh14} introduced a Godunov-type scheme for the shallow water equations based on a Riemann solver. For some other numerical methods on resonance systems, see \cite{Gallouet,Saurel17}.

 \begin{figure}[h]
\centering
\includegraphics[width=0.56\textwidth]{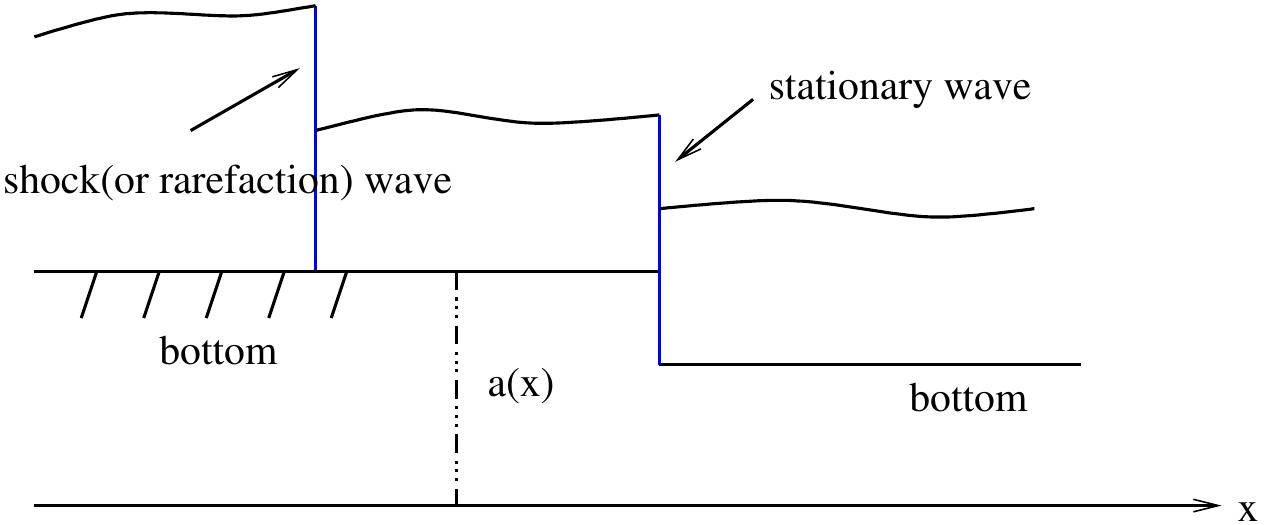}
\caption*{Fig. 1.2. The interaction of water waves.}
\end{figure}

The water wave types are complicated in the shallow water system since different waves may break or merge into other waves, e.g., roll waves and solitary waves. Thus it is an important topic to classify the interaction results of shallow water system. One of the basic questions is the interaction of elementary waves, see Fig. 1.2. Since the interaction results of elementary waves apart from the stationary wave for gas dynamics have already been  obtained by Chang and Hsiao (\cite{ChangHsiao3}) as well as Smoller \cite{Smoller}, one can directly apply the results to the shallow water system without extra difficulties. Our objective here is to discuss the rarefaction waves or shock waves interacting with the stationary wave. Recently, Sheng and Zhang \cite{Sheng} have investigated the interaction of elementary waves in a variable cross section area, especially for rarefaction wave or shock wave interacting with a stationary wave. Here our main goal is to extend their results to the shallow water equations \eqref{1.1}. On the one hand, by using the characteristic analysis methods, we prove that when a forward rarefaction wave interacts with a stationary wave, it will transmit a forward rarefaction wave. At the same time, a backward wave will transmit or reflect as well.
 Furthermore, we consider a more complicated case in which a compressible wave is transmitted. The results in large time scales are well investigated.  On the other hand,  when a forward shock wave interacts with a stationary wave, it will penetrate the stationary wave, and either transmit or reflect a backward wave. We believe our results can apply to other systems with similar structures.

The rest of this paper is organized as follows. In Section \ref{sec2}, we perform characteristic analysis for the shallow water equations, and discuss the properties of elementary waves. In Section \ref{sec3},  the interactions of rarefaction wave and shock wave with stationary wave are discussed. By using the  characteristic analysis methods, we are able to describe the evolution of waves during the interaction process and for large time scale in each case.

\section{Preliminaries}\label{sec2}
\subsection{Characteristic analysis and elementary waves}
Denote $U=(u,h,a)$, system \eqref{1.1} can be rewritten, when considering a smooth solution, as
\begin{equation}\label{2.1}
A(U)\partial_t U+B(U)\partial_x U=0,
\end{equation}
where
$$
A(U)=\left( \begin{array}{ccc}
0 & 0 & h\\
0 & 1 & 0 \\
1 & 0 & 0 \end{array} \right),
\quad
B(U)=\left( \begin{array}{ccc}
-gh & gh & hu \\
0 & u & h \\
0 & 0 & 0 \end{array} \right).
$$
\eqref{2.1} has three eigenvalues
\begin{equation}\label{2.2}
\lambda_1=u-c,\quad  \lambda_2=u+c,\quad \lambda_3=0,
\end{equation}
where $c$ is the celerity denoted by (\cite{Toro})
\begin{equation}\label{2.2'}
c=\sqrt{gh}.
\end{equation}
The corresponding right eigenvectors are:
\begin{equation}\label{2.3}
\vec{r}_1=\left(h,-c,0\right)^T,\quad \vec{r}_2=\left(h,c,0\right)^T,\quad \vec{r}_3=\left(c^2,-gu,u^2-c^2\right)^T.
\end{equation}
The third characteristic family is linearly degenerate, while the first and the second characteristics families are genuinely nonlinear:
\begin{equation}\label{2.4}
-\triangledown \lambda_1(U)\cdot r_1(U)=\triangledown \lambda_2(U)\cdot r_2(U)=\frac32c>0.
\end{equation}
One may notice that the first and the second characteristics may coincide with the third one, so the system is not strictly hyperbolic. More precisely, setting
\begin{equation}\label{2.5}
\Gamma_{\pm}:u=\pm c,
\end{equation}
one has
\begin{equation}\label{2.6}
\lambda_1=\lambda_3\quad {\rm on} \quad \Gamma_{+},\qquad
\lambda_2=\lambda_3\quad {\rm on} \quad  \Gamma_{-}.
\end{equation}
In the $(u,h,a)$ space, $\Gamma_{\pm}$ separate the half-plane $h>0$ into three parts. For convenience, we will view them as $D_{1}, D_2$ and $D_3$:
\begin{equation}\label{2.7'}
D_1  =\Big\{(u,h,a)\big|u<-c~\Big\},  \quad
D_2  =\Big\{(u,h,a)\big||u|<c~\Big\},\quad
D_3  =\Big\{(u,h,a)\big|u>c~\Big\},
\end{equation}
To further simplify the descriptions, we call $D_2$ as $\emph {subcritical area} $($|u|<c$), and $D_1, D_3$ as $\emph {supercritical area}$ ($|u|>c$).
It is also convenient to set \\
\begin{equation}\label{2.8}
D_2^{+}:=\Big\{(u,h,a)\in D_2, u\geq 0\Big\},\quad D_2^{-}:=\Big\{(u,h,a)\in D_2, u< 0\Big\}.
\end{equation}
In each of the regions, the system is strictly hyperbolic and
\begin{equation}\label{2.9}
\begin{array}{llll}
&\lambda_1<\lambda_2<\lambda_3, &{\rm in} & D_1,\\
&\lambda_1<\lambda_3<\lambda_2, &{\rm in} & D_2,\\
&\lambda_3<\lambda_1<\lambda_2, &{\rm in} & D_3.
\end{array}
\end{equation}

\subsection{The rarefaction wave}
First, we look for self-similar solutions. The Riemann invariants of each characteristic are calculated as :
\begin{equation}\label{2.10}
\left\{
\begin{array}{lll}
&\lambda_1=u-\sqrt{gh}: &\big\{a, u+2\sqrt{gh}\big\},  \\[5pt]
&\lambda_2=u+\sqrt{gh}: &\big\{a, u-2\sqrt{gh}\big\},   \\[5pt]
&\lambda_3=0: &\displaystyle \big\{hu, \frac{u^2}2+g(h+a)\big\}.
\end{array}
\right.
\end{equation}
For the rarefaction waves, the variable $a(x)$ remains constant,  system \eqref{1.1} degenerates to the shallow water equations without a bottom step
\begin{equation}\label{2.11}
\left\{
\begin{array}{l}
h_t+(h u)_{x}=0,\\[5pt]
\displaystyle (h u)_t+(hu^2+\frac{gh^2}2)_x=0.\\
\end{array}
\right.
\end{equation}
For a given left-hand state $U_0=(u_0,h_0,a_0)$, the right-hand states $U$ that can be connected by 1-wave and 2-wave rarefaction curves are determined by
\begin{equation}\label{2.12}
\left\{
\begin{array}{l}
\displaystyle \overleftarrow{R}_1(U,U_0):~ u=u_0-2\sqrt{g}\left(\sqrt{h}-\sqrt{h_0}\right),\quad h<h_0,\\[8pt]
\displaystyle \overrightarrow{R}_2(U,U_0):~ u=u_0+2\sqrt{g}\left(\sqrt{h}-\sqrt{h_0}\right),\quad h>h_0.
\end{array}
\right.
\end{equation}
\subsection{The stationary wave}
The Rankine-Hugoniot relation associated with the third equation of \eqref{1.1} is that
\[-\sigma[a]=0,\]
where $[a]:~=a_1-a_0$ is the jump of the variable $a$. This leads to the following conclusions:\\
1)~~$\sigma=0:$~the shock speed vanishes, here assume $[a]\not=0$ and called stationary contact discontinuity; \\
2)~~$[a]=0:$~the bottom level~$a$~remains constant across the non-zero speed shocks.\\
Across the stationary contact discontinuity, the Riemann invariants remain constant. From the third equation of \eqref{2.10}, one has 
\begin{equation}\label{2.13}
\left\{
\begin{array}{l}
[hu]=0,\\[5pt]
\displaystyle [ \frac{u^2}2+g(h+a)]=0.
\end{array}
\right.
\end{equation}
\eqref{2.13} determines the stationary contact curve which is parameterized as
\begin{equation}\label{2.14}
\left\{
\begin{array}{l}
\displaystyle u=\frac{h_0u_0}{h},\\[5pt]
\displaystyle a=a_0+\frac{u_0^2-u^2}{2g}+h_0-h,
\end{array}
\right.
\end{equation}
where $(u_0,h_0,a_0)$ is the given left-hand state. Moreover, we have the following lemma.
\begin{lem}\label{lem1}
Given the left-hand state $U_0=(u_0,h_0,a_0)$, \eqref{2.14} has at most two solutions $U_*=(u_*, h_*, a)$ and $U^{*}=(u^{*}, h^{*}, a)$ for any $a>0$, if and only if $a<\displaystyle a_{\rm max}(U_0)$, where
$$
a_{\rm max}(U_0)=a_0+h_0+\frac{u_0^2}{2g}-\frac{3}{2g^{1/3}}(h_0u_0)^{2/3}.
$$
More precisely,\\
1) If $a>a_{\rm max}(U_0)$, \eqref{2.14} has no solution, so there is no stationary contact.\\
2) If $a< a_{\rm max}(U_0)$, there are two points $U_*, U^{*}$ satisfying \eqref{2.14}, which can be connected with $U_0$ by a stationary contact.\\
3) If $a=a_{\rm max}(U_0)$, $U_*$ and $U^{*}$ coincide.
\end{lem}
%\begin{proof}
%We substitute the first equation of \eqref{2.8} into the second equation and get
%\begin{equation}\label{2.9}
%a_0-a+\frac1{2g}\left(u_0^2-\left(\frac{h_0u_0}{h}\right)^2\right)+h_0-h=0.
%\end{equation}
%Let
%$$
%\displaystyle \Phi(h)=\Phi(U_0,a;h)=h^3+\left(a-a_0-h_0-\frac{u_0^2}{2g}\right)h^2+\frac{h_0^2u_0^2}{2g},
%$$
%we have
%\begin{equation}\label{2.10}
%\frac{{\rm d}\Phi(h) }{{\rm d}h}=3h^2+2h\left(a-a_0-h_0-\frac{u_0^2}{2g}\right),
%\end{equation}
%as $h>0$, $\Phi(h)$ reaches its minimum value at
%$$
%\bar{h}_m(U_0)=\frac{u_0^2+2g(a_0+h_0-a)}{3g}.
%$$
%So $\Phi(h)=0$ admits a solution when $\Phi(\bar{h}_m)\leq 0$.
%It follows that
%$$
%a\leq a_{\rm max}(U_0)=a_0+h_0+\frac{u_0^2}{2g}-\frac{3}{2g^{1/3}}(h_0u_0)^{2/3}.
%$$
%If $a< a_{\rm max}(U_0)$, then there exist exactly two values $h_\ast(U_0)<\bar{h}_m(U_0)<h^{\ast}(U_0)$ which satisfy
%\begin{equation}\label{2.11}
%\Phi(U_0,a;h)=\Phi(U_0,a;h)=0,
%\end{equation}
%the two values coincide iff $a=a_{\rm max}(U_0)$. Thus we complete the proof of the lemma 2.1.
%\end{proof}
The proof of lemma \ref{lem1} is straightforward, see \cite{LeflochThanh13} for details. Moreover, across the stationary contact discontinuity denoted by $S_0(U,U_0)$, the states $U_*=(u_*,h_*,a)$ and $U^*=(u^{*},h^{*},a)$ have the following properties
$$
S_0(U,U_0)=\left\{\begin{array}{lll}
S_0(U_*,U_0), &|u_*|>\sqrt{gh_*},\\
S_0(U^*,U_0), &|u^*|<\sqrt{gh^*},
\end{array}\right.
{\rm more~ precisely},
S_0(U,U_0):~\left\{\begin{array}{lll}
 U_* \in\left\{
 \begin{array}{lll}
 D_1,&u_0<0,\\
 D_3,&u_0>0,
 \end{array}\right.\\
U^{*} \in D_2.
\end{array}\right.
$$
For the proof of the properties, we refer to \cite{LeflochThanh14}.

As shown in \cite{LeflochThanh13,LeflochThanh14}, the Riemann problem for \eqref{1.1} may admit up to a one-parameter family of solutions. This phenomena can be avoided by requiring Riemann solutions to satisfy an admissibility criterion: monotone condition on the component $a$. Following \cite{LeflochThanh13,LeflochThanh14,Sheng}, one can impose the global entropy condition on stationary contact discontinuity of \eqref{1.1}.\\
{\bf Global entropy condition.} Along the stationary curve $S_0(U,U_0)$ in the $(u, h)$-plane, $a(h)$ obtained from \eqref{2.14} is a monotone function of $h$.

Under the global entropy condition, the stationary contact discontinuity can be called as stationary wave.  LeFloch and Thanh further \cite{LeflochThanh12} proved the following results.
\begin{lem}\label{lem2}
Global entropy condition is equivalent to the statement that any
stationary wave has to remain in the closure of only one domain $D_i, i=1,2,3$. More precisely, \\
1) If $U_0\in D_1\cup D_3$, $h=h_{*}(U_0)$ is chosen as the solution.\\
2) If $U_0 \in D_2$, $h=h^{*}(U_0)$ is chosen as the solution.
\end{lem}
From lemma \ref{lem2}, we next investigate the properties of the stationary curve in the $(u,h)$-plane.
\begin{lem}\label{lem3}
The stationary wave can be viewed as parametrized curves $S_0(U(a);U_0)$ depending only on the bottom level $a$ in $(u,h)$ plane, and they have the following properties:

1) $S_0(U(a);U_0)$ is a convex curve which is strictly increasing (decreasing) in $u$ if $u<0~(>0)$.

2) The relation of $U_0$ and $U(a)$ is shown as\\
$$
U_0\in D_1:~\left\{\begin{array}{lll}
u>u_0, h>h_0 &a>a_0,\\
u<u_0, h<h_0 &a<a_0,
\end{array}\right.,~~
U_0\in D_3:~\left\{\begin{array}{lll}
u<u_0, h>h_0 &a>a_0,\\
u>u_0, h<h_0 &a<a_0,
\end{array}\right.
$$
$$
U_0\in D_2^{-}:~\left\{\begin{array}{lll}
u<u_0, h<h_0 &a>a_0,\\
u>u_0, h>h_0 &a<a_0,
\end{array}\right.,~~
U_0\in D_2^{+}:~\left\{\begin{array}{lll}
u>u_0, h<h_0 &a>a_0,\\
u<u_0, h>h_0 &a<a_0.
\end{array}\right.
$$
\end{lem}

\begin{figure}[htbp]
\centering
\includegraphics[width=0.66\textwidth]{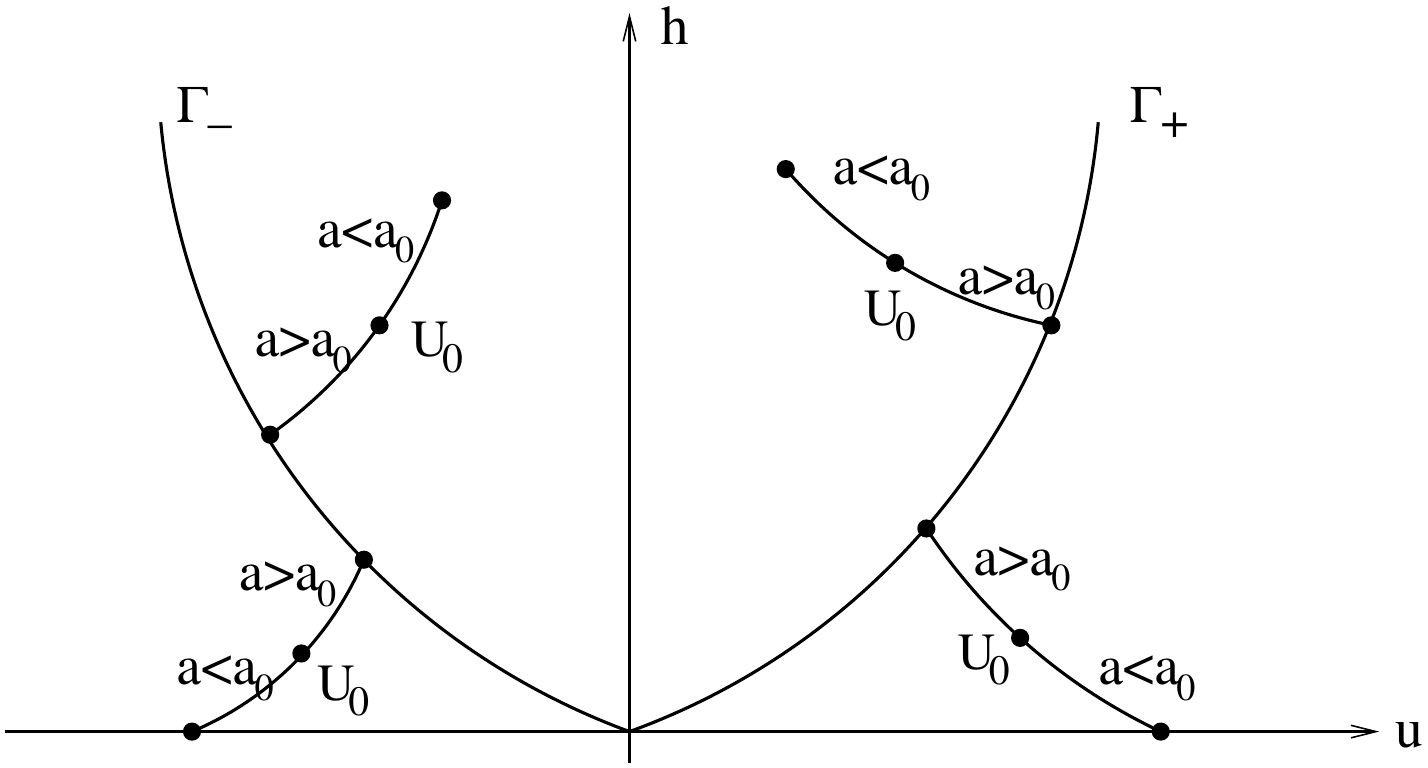}
\caption*{Fig. 2.1. The stationary curve $S_0(U,U_0)$ in $(u,h)$ plane given $U_0$.}
\end{figure}

\begin{proof}
Differentiating the two equations of \eqref{2.14}, one gets
\begin{equation}\label{2.15}
\left\{
\begin{array}{ll}
\displaystyle h{\rm d}u+u{\rm d}h=0,\\
u{\rm d}u+g({\rm d}h+{\rm d}a)=0,
\end{array}
\right.
\end{equation}
which leads to
\begin{equation}\label{2.16}
\begin{array}{lll}
\displaystyle \frac{{\rm d} h}{{\rm d} u}=-\frac{h}{u},\quad \frac{{\rm d}a}{{\rm d}h}=\frac{u^2-gh}{gh},\quad \frac{{\rm d}a}{{\rm d}u}=-\frac{u^2-gh}{gu}.
\end{array}
\end{equation}
Statement 1) is derived from the first formula of \eqref{2.16}, and statement 2) follows the second and the third formulas of \eqref{2.16}. See Fig. 2.1. Thus we prove lemma \ref{lem3}.
\end{proof}

Under the following transformation
$$
x \rightarrow -x,\quad u\rightarrow -u,
$$
a right-hand state $U=(h,u,a)$ is transformed into a left-hand state $U^{'}=(h,-u,a)$. Without loss of generality, we assume that the bottom level is decreasing from now on, i.e., $a_0>a_1$.

\subsection{The shock wave}
For the non-zero speed shocks, the Rankine-Hugoniot relations corresponding to \eqref{2.11} are 
\begin{equation}\label{2.17}
\left\{\begin{array}{ll}
-\sigma[h]+[h u]=0,\\
\displaystyle -\sigma[hu]+[hu^2+\frac{gh^2}2]=0,
\end{array}\right.
\end{equation}
which  is equivalent to
\begin{equation}\label{2.18}
\sigma_i(U, U_0)=u_0\mp \sqrt{\frac{g}2(h+h_0)\frac{h}{h_0}},\quad i=1,2.
\end{equation}
the 1-and 2-families of shock waves connecting a given left-hand state $U_0$ to the right-hand state $U$ is constrained by the Hugoniot set
\begin{equation}\label{2.19}
(u-u_0)^{2}=\frac{g}2(h-h_0)^2\left(\frac1{h}+\frac1{h_0}\right).
\end{equation}
A shock wave should also satisfy the Lax shock conditions
\begin{equation}\label{2.20}
\lambda_i(U)<\sigma_i(U,U_0)<\lambda_i(U_0),\quad i=1,2.
\end{equation}
Thus the 1- and 2-shock waves 
$\overleftarrow{S}_1(U_0)$  and $\overrightarrow{S}_2(U_0)$ consisting of all right-hand states $U$ are determined by
\begin{equation}\label{2.21}
\left\{
\begin{array}{l}
\displaystyle \overleftarrow{S}_1(U_0): u=u_0-\sqrt{\frac{g}2(h-h_0)^2\left(\frac1{h}+\frac1{h_0}\right)},\quad h>h_0,\\
\displaystyle \overrightarrow{S}_2(U_0): u=u_0-\sqrt{\frac{g}2(h-h_0)^2\left(\frac1{h}+\frac1{h_0}\right)},\quad h<h_0.
\end{array}
\right.
\end{equation}
The 1- and 2-shock wave speeds ~$\sigma_i(U,U_0) (i=1, 2)$ may change their signs along the shock curves in the $(u, h)$ plane, more precisely,
\begin{equation}\label{2.22}
\sigma_1(U,U_0)~\left\{\begin{array}{lll}
~<0, &U_0\in D_1\cup D_2,\\
\left.\begin{array}{ll}
<0, &\tilde{h}_0<h,\\
=0, &h=\tilde{h}_0,\\
>0, &h_0<h<\tilde{h}_0,
\end{array}\right\} &U_0\in D_3,
\end{array}\right.
\end{equation}
and
\begin{equation}\label{2.23}
\sigma_2(U,U_0)~\left\{\begin{array}{lll}
~>0, &U_0\in D_2\cup D_3,\\
\left.\begin{array}{ll}
>0, &\bar{h}_0<h,\\
=0, &h=\bar{h}_0,\\
<0, &h_0<h<\bar{h}_0,
\end{array}\right\} &U_0\in D_1,
\end{array}\right.
\end{equation}
where $\widetilde{U}_0=(\tilde{u}_0, \tilde{h}_0,a_0)\in D_2\cap\{u>0\}, \overline{U}_0=( \bar{u}_0, \bar{h}_0,a_0)\in D_2\cap\{u<0\}$.

For convenience, let us define the backward and forward wave curves
\begin{align*}
W_1(h;U_0) & =\left\{\begin{array}{ll}
\overleftarrow{R_1}(h;U_0), &h<h_0,\\
\overleftarrow{S_1}(h;U_0), &h>h_0,
\end{array}\right.\\
W_2(h;U_0) & =\left\{\begin{array}{ll}
\overrightarrow{R_2}(h;U_0), &h>h_0,\\
\overrightarrow{S_2}(h;U_0), &h<h_0,
\end{array}\right.
\end{align*}
and stationary wave
$$
W_3(h;U_0)=S_0(h;U_0),~~ h=h_*~~{\rm or}~~h^*.
$$
A straight calculation shows that the wave curve $W_1(h;U_0)$ is strictly decreasing and convex in the $(u, h)$ plane, while the wave curve $W_2(h;U_0)$ is strictly increasing and convex.

The elementary waves of system \eqref{1.1} consist of rarefaction waves, shock waves, and stationary wave, which are denoted by $W_i(h;U_0)~(i=1,2,3)$ briefly.

According to the Riemann solutions (\cite{LeflochThanh13}) of \eqref{1.1} in $D_1, D_2$, and $D_3$ , we will consider the interaction results of the elementary waves based on the division of the $(u, h)$ plane. Here we only consider the interactions of rarefaction wave or shock wave with the stationary wave.

\section{The interactions of rarefaction wave or shock wave with the stationary wave}\label{sec3}

To study the interactions of rarefaction wave or shock wave with the stationary wave, we consider the initial value problem \eqref{1.1} with
\begin{equation}\label{3.1}
(u,h,a)\Big|_{t=0}=
\left\{\begin{array}{lll}
U_-=(u_-,h_-,a_0), &x<x_1,\\
U_m=(u_m,h_m,a_0), &x_1<x<x_2,\\
U_+=(u_+,h_+,a_1), &x>x_2.
\end{array}\right.
\end{equation}
First, we investigate the interaction of a rarefaction wave with a stationary wave. For definiteness, we make the assumption that $U_-, U_m$ and $U_+$ are in the first quadrant of $(u,h)$ plane.
\subsection{Interaction of a rarefaction wave with a stationary wave}
 In this case, $U_m \in \overrightarrow{R}_2(h;U_-)$, $U_+ \in W_3(h;U_m)$, $a_1< a_{max}(U_m)$,
while $\displaystyle a_{max}(U_m)=a_0+h_m+\frac{u_m^2}{2g}-\frac3{2g^{1/3}}\left(h_mu_m\right)^{2/3}$. A straight calculation shows that $a<a_{max}(U_0)$ for $h_-\leq h_0\leq h_m$, which means the stationary wave will always exist in the interaction process. We have that
\begin{equation}\label{3.2}
\left\{\begin{array}{ll}
\overrightarrow{R}_2(U,U_-)=\overrightarrow{R}_2(U_m,U):~ u=u_-+2\sqrt{g}(\sqrt{h}-\sqrt{h_-}),\qquad ~~h_-\leq h\leq h_m,\\[5pt]
W_3(U_1;U_0):~\left\{\begin{array}{lll}
h_0 u_0 = h_1u_1, &U_0=(u_0,h_0,a_0)\in \overrightarrow{R}_2(U_m,U),\\
\displaystyle \frac{u_0^2}2+g(h_0+a_0)=\frac{u_1^2}2+g(h_1+a_1), &a_1<a_{max}.
\end{array}
\right.
\end{array}\right.
\end{equation}
First, we use two curves $\Gamma_+$ and $\Gamma_1$ to divide the first quadrant of $(u,h)$ plane into three parts {\bf I, I\!I} and {\bf I\!I\!I}, where $\Gamma_1: u=2c$, namely (see Fig. 3.1)
\begin{equation}\label{3.3}
\begin{array}{cc}
{\rm I}=\big\{(u,\rho)\big|u>2c\big\}, \quad
{\rm I\!I}= \big\{c<u<2c\big\}, \\[5pt]
{\rm I\!I\!I}= \big\{0<u<c\big\}. \quad
\end{array}
\end{equation}

%TeXCAD Picture [three_parts.pic]. Options:
%\grade{\on}
%\emlines{\off}
%\epic{\off}
%\beziermacro{\on}
%\reduce{\on}
%\snapping{\off}
%\quality{8.00}
%\graddiff{0.01}
%\snapasp{1}
%\zoom{6.7272}
\unitlength 0.75mm % = 2.85pt
\linethickness{0.4pt}
\ifx\plotpoint\undefined\newsavebox{\plotpoint}\fi % GNUPLOT compatibility
\begin{picture}(108.37,85.24)(-26,5)
\footnotesize
%\vector(18.08,22.07)(106.27,22.18)
\put(106.27,22.18){\vector(1,0){.07}}\put(18.08,22.07){\line(1,0){88.19}}
%\end
\put(23.65,15.45){\vector(0,1){69.79}}
\put(26,11.05){Fig. 3.1. The three parts in $(u, h)$ plane: ${\bf I, I\!I, I\!I\!I}$.}
\qbezier(23.69,22.01)(54.09,35.44)(59.75,83.17)
\qbezier(23.69,22.01)(81.67,29.83)(86.97,83.79)
\put(108.37,21.85){$u$}
\put(21.11,84.58){$h$}
\put(55,81){$\Gamma_+$}
\put(82,81){$\Gamma_1$}
\put(66.89,57.97){${\bf I\!I}$}
\put(38.5,57.82){${\bf I\!I\!I}$}
\put(96.33,57.38){${\bf I}$}
\end{picture}

\begin{lem}\label{lem3.1}
There hold the following properties (see Fig. 3.2):

1) If $U_+$ $\big(U_m\big)$ is on the left of $\Gamma_1$, then $\overrightarrow{R}_2(U_+,U)$ $\big(\overrightarrow{R}_2(U_m,U)\big)$ intersects with $h-$ axis at  $\bar{h}_+$ $(\bar{h}_m)$.

2) If $U_+$ $\big(U_m\big)$ is on the right of $\Gamma_1$, then $\overrightarrow{R}_2(U_+,U)$ $\big(\overrightarrow{R}_2(U_m,U)\big)$ intersects with $u-$ axis at  $\bar{u}_+$ $(\bar{u}_m)$.
\end{lem}

%TeXCAD Picture [three_parts''.pic]. Options:
%\grade{\on}
%\emlines{\off}
%\epic{\off}
%\beziermacro{\on}
%\reduce{\on}
%\snapping{\off}
%\quality{8.00}
%\graddiff{0.01}
%\snapasp{1}
%\zoom{16.0001}
\unitlength 0.8mm % = 2.85pt
\linethickness{0.4pt}
\ifx\plotpoint\undefined\newsavebox{\plotpoint}\fi % GNUPLOT compatibility
\begin{picture}(108.37,85.24)(-25,5)
\footnotesize
\put(23.65,15.45){\vector(0,1){69.79}}
\put(31.55,11.05){Fig. 3.2. The properties of the curve $\Gamma_1$.}
\put(108.37,21.85){$u$}
\put(20,84){$h$}
\qbezier(23.64,22.15)(69.35,36.72)(73.14,83.39)
\put(68.5,82.6){$\Gamma_1$}
%\dashline{1}(81.16,61.1)(91.87,53.07)
\multiput(81.09,61.03)(.044625,-.033458){16}{\line(1,0){.044625}}
\multiput(82.52,59.96)(.044625,-.033458){16}{\line(1,0){.044625}}
\multiput(83.95,58.89)(.044625,-.033458){16}{\line(1,0){.044625}}
\multiput(85.37,57.82)(.044625,-.033458){16}{\line(1,0){.044625}}
\multiput(86.8,56.75)(.044625,-.033458){16}{\line(1,0){.044625}}
\multiput(88.23,55.68)(.044625,-.033458){16}{\line(1,0){.044625}}
\multiput(89.66,54.61)(.044625,-.033458){16}{\line(1,0){.044625}}
\multiput(91.09,53.54)(.044625,-.033458){16}{\line(1,0){.044625}}
%\end
\put(82.05,63.18){$U_m$}
\put(94.1,54.11){$U_+$}
\put(18.14,22){\vector(1,0){88.3}}
%\dashline{1}(49.95,70.76)(55,60.5)
\multiput(49.88,70.69)(.03237,-.06577){12}{\line(0,-1){.06577}}
\multiput(50.66,69.11)(.03237,-.06577){12}{\line(0,-1){.06577}}
\multiput(51.43,67.53)(.03237,-.06577){12}{\line(0,-1){.06577}}
\multiput(52.21,65.95)(.03237,-.06577){12}{\line(0,-1){.06577}}
\multiput(52.99,64.38)(.03237,-.06577){12}{\line(0,-1){.06577}}
\multiput(53.76,62.8)(.03237,-.06577){12}{\line(0,-1){.06577}}
\multiput(54.54,61.22)(.03237,-.06577){12}{\line(0,-1){.06577}}
%\end
\qbezier(54.85,60.5)(40.43,47.79)(23.64,45.49)
\qbezier(49.95,70.61)(40.21,59.01)(23.64,55.74)
\put(51.73,71.5){$U_m$}
\put(57.23,60.35){$U_+$}
{\tiny
\put(69,44.36){$\overrightarrow{R}_2(U_m,U)$}
\put(85,35){$\overrightarrow{R}_2(U_+,U)$}
\put(36,62.54){$\overrightarrow{R}_2(U_m,U)$}
\put(40,49.82){$\overrightarrow{R}_2(U_+,U)$}}
\put(18,55.6){$\bar{h}_m$}
\put(18,45.41){$\bar{h}_+$}
\qbezier(81.19,61.06)(66.97,25.34)(56.87,22)
\qbezier(91.81,53.06)(83.66,28.56)(74.87,22.06)
\put(75,18){$\bar{u}_+$}
\put(57.06,18){$\bar{u}_m$}
\end{picture}

As shown in \cite{ChangHsiao3}, here we use phase plane analysis method to discuss the interaction results. Recall that $S_0(U_1,U_0)$ is the right-hand state $U_1=(u_1,h_1,a_1)$ starting from $U_0=(u_0,h_0,a_0)$ by a stationary wave. If we set $U_0\in\overrightarrow{R}_2(U_-,U_m), 0<h_0<h_m$, then $S_0(U_1,U_0)$ starts from $U_1=U_+$ as $U_0=U_m$.  The key point here is to compare the relative positions between the curves $S_0(U_1,U_0)$ and $\overrightarrow{R}_2(U_+,U)$. Since both curves can be written as parameterized functions of the variable $h_1$, we can compare them directly and have the following lemma. 

%\begin{equation}\label{3.3}
%\left\{
%\begin{array}{ll}
%S_0(h_1;U_0):=\big\{u_1(h_1;U_0)\big|(u_1,h_1,a_1)\in S_0(U,U_0)\big \},\\[5pt]
%U_0\in\overrightarrow{R}_2(U_-;U_m),~~~~~0<h_0<h_m.
%\end{array}
%\right.
%\end{equation}
%then $S_0(h_1;U_0)$ is a curve parameterized by $h_1$. From \eqref{3.2}, it is obviously that  $S_0(h_1;U_0)$ passes through $U_+$. The key point here is to compare the relative position between $S_0(h_1;U_0)$ and $\overrightarrow{R}_2(U_+,U)$, which gives us a clue to analyze the interaction results. 

\begin{lem}\label{lem3.2}
The relative positions of the curves $S_0(U_1,U_0)$ and $\overrightarrow{R}_2(U_+,U_1)$, where $U_0=U_m$, $U_1=U_+$, are shown as follows.

1) When $U_+\in {\rm I}$, $U_m\in {\rm I}$, $S_0(U_1,U_0)$ is below (or on the right side of) the curve $\overrightarrow{R}_2(U_+,U)$ in the half plane $h>0$ (see Fig .3.3);

2) When $U_+\in {\rm I}\cup {\rm I\!I}$, $U_m\in {\rm I\!I}$, $S_0(U_1,U_0)$ is below the curve $\overrightarrow{R}_2(U_+,U)$ in the supercritical area, and is above $\overrightarrow{R}_2(U_m,U)$ in the subcritical area(see Fig .3.4-3.5);

3) When $U_+\in {\rm I\!I\!I}$, $U_m\in {\rm I\!I\!I}$, $S_0(U_1,U_0)$ is below the curve $\overrightarrow{R}_2(U_+,U)$ (see Fig .3.6).

\end{lem}

%TeXCAD Picture [traffic flow 2.pic]. Options:
%\grade{\on}
%\emlines{\off}
%\epic{\off}
%\beziermacro{\on}
%\reduce{\on}
%\snapping{\off}
%\quality{8.00}
%\graddiff{0.01}
%\snapasp{1}
%\zoom{6.7272}
\unitlength 0.7mm % = 2.85pt
\linethickness{0.4pt}
\ifx\plotpoint\undefined\newsavebox{\plotpoint}\fi % GNUPLOT compatibility
\begin{picture}(108.75,83.75)(32,6)
\footnotesize
\put(26.5,22.75){\vector(1,0){80}}
\put(31,19.75){\vector(0,1){64}}
\put(108.75,22.5){$u$}
\qbezier(31,22.75)(52.88,36.25)(54.25,78.75)
\qbezier(31,22.75)(74.38,31.75)(77.25,78.75)
%\dashline{1}(79.75,64.25)(101,59)
\multiput(79.68,64.18)(.13199,-.03261){7}{\line(1,0){.13199}}
\multiput(81.53,63.72)(.13199,-.03261){7}{\line(1,0){.13199}}
\multiput(83.38,63.27)(.13199,-.03261){7}{\line(1,0){.13199}}
\multiput(85.22,62.81)(.13199,-.03261){7}{\line(1,0){.13199}}
\multiput(87.07,62.35)(.13199,-.03261){7}{\line(1,0){.13199}}
\multiput(88.92,61.9)(.13199,-.03261){7}{\line(1,0){.13199}}
\multiput(90.77,61.44)(.13199,-.03261){7}{\line(1,0){.13199}}
\multiput(92.61,60.98)(.13199,-.03261){7}{\line(1,0){.13199}}
\multiput(94.46,60.53)(.13199,-.03261){7}{\line(1,0){.13199}}
\multiput(96.31,60.07)(.13199,-.03261){7}{\line(1,0){.13199}}
\multiput(98.16,59.61)(.13199,-.03261){7}{\line(1,0){.13199}}
\multiput(100.01,59.16)(.13199,-.03261){7}{\line(1,0){.13199}}
%\end
\qbezier(80,64.25)(74.13,38.5)(55.75,22.75)
\qbezier(100.5,59)(96.38,29.63)(89.75,22.75)
\put(80,64){\circle*{1}}
\put(100.5,58.75){\circle*{1}}
%\dashline{1}(76,50.75)(97.75,42.25)
\multiput(75.93,50.68)(.08239,-.0322){11}{\line(1,0){.08239}}
\multiput(77.74,49.97)(.08239,-.0322){11}{\line(1,0){.08239}}
\multiput(79.55,49.26)(.08239,-.0322){11}{\line(1,0){.08239}}
\multiput(81.37,48.55)(.08239,-.0322){11}{\line(1,0){.08239}}
\multiput(83.18,47.85)(.08239,-.0322){11}{\line(1,0){.08239}}
\multiput(84.99,47.14)(.08239,-.0322){11}{\line(1,0){.08239}}
\multiput(86.8,46.43)(.08239,-.0322){11}{\line(1,0){.08239}}
\multiput(88.62,45.72)(.08239,-.0322){11}{\line(1,0){.08239}}
\multiput(90.43,45.01)(.08239,-.0322){11}{\line(1,0){.08239}}
\multiput(92.24,44.3)(.08239,-.0322){11}{\line(1,0){.08239}}
\multiput(94.05,43.6)(.08239,-.0322){11}{\line(1,0){.08239}}
\multiput(95.87,42.89)(.08239,-.0322){11}{\line(1,0){.08239}}
%\end
\put(76,50.75){\circle*{1}}
\put(97.75,42){\circle*{1}}
\put(27,81){$h$}
\put(81.5,65.75){$U_m$}
\put(102,60.75){$U_+$}
\put(79,50){$U_0$}
\put(99,40.5){$U_1$}
\put(36.5,10.25){Fig. 3.3. $U_+\in {\rm I}$, $U_m\in {\rm I}$.}
{\tiny
\put(90,29.5){$S_0(U_1,U_0)$}
\put(61,34){$\overrightarrow{R}_2(U_m,U)$}
\put(80,37.25){$\overrightarrow{R}_2(U_+,U)$}}
\qbezier(100.41,58.69)(90.91,30.85)(77.69,22.8)
\put(87.26,71.95){${\rm I}$}
\put(66.3,71.35){${\rm I\!I}$}
\put(44.74,71.35){${\rm I\!I\!I}$}
\put(90.25,18){$\tilde{u}$}
\put(55.6,18){$\bar{u}_m$}
\put(77.74,18){$\bar{u}_+$}
\put(56.78,79.23){$\Gamma_+$}
\put(79.68,79.08){$\Gamma_1$}
\end{picture}
%TeXCAD Picture [traffic flow 3.pic]. Options:
%\grade{\on}
%\emlines{\off}
%\epic{\off}
%\beziermacro{\on}
%\reduce{\on}
%\snapping{\off}
%\quality{8.00}
%\graddiff{0.01}
%\snapasp{1}
%\zoom{6.7272}
\unitlength 0.7mm % = 2.85pt
\linethickness{0.4pt}
\ifx\plotpoint\undefined\newsavebox{\plotpoint}\fi % GNUPLOT compatibility
\begin{picture}(50.75,50.75)(44,6)
\footnotesize
\put(26.5,22.75){\vector(1,0){80}}
\put(31,19.75){\vector(0,1){64}}
\put(108.75,22.5){$u$}
\qbezier(31,22.75)(52.88,36.25)(54.25,78.75)
\put(27,80){$h$}
\put(36.5,10.25){Fig. 3.4. $U_+\in {\rm I}$, $U_m\in {\rm I\!I}$.}
%\dashline{1}(72.92,69.21)(91.13,61.78)
\multiput(72.85,69.14)(.07883,-.03216){11}{\line(1,0){.07883}}
\multiput(74.58,68.43)(.07883,-.03216){11}{\line(1,0){.07883}}
\multiput(76.32,67.72)(.07883,-.03216){11}{\line(1,0){.07883}}
\multiput(78.05,67.02)(.07883,-.03216){11}{\line(1,0){.07883}}
\multiput(79.79,66.31)(.07883,-.03216){11}{\line(1,0){.07883}}
\multiput(81.52,65.6)(.07883,-.03216){11}{\line(1,0){.07883}}
\multiput(83.26,64.89)(.07883,-.03216){11}{\line(1,0){.07883}}
\multiput(84.99,64.19)(.07883,-.03216){11}{\line(1,0){.07883}}
\multiput(86.72,63.48)(.07883,-.03216){11}{\line(1,0){.07883}}
\multiput(88.46,62.77)(.07883,-.03216){11}{\line(1,0){.07883}}
\multiput(90.19,62.06)(.07883,-.03216){11}{\line(1,0){.07883}}
%\end
\qbezier(91.13,61.78)(74.25,32.7)(55.24,22.72)
\qbezier(72.92,69.12)(59.22,46.98)(31.02,36.68)
%\dashline{1}(48.35,44.9)(77.16,30.76)
\multiput(48.28,44.83)(.067156,-.03296){13}{\line(1,0){.067156}}
\multiput(50.03,43.97)(.067156,-.03296){13}{\line(1,0){.067156}}
\multiput(51.77,43.12)(.067156,-.03296){13}{\line(1,0){.067156}}
\multiput(53.52,42.26)(.067156,-.03296){13}{\line(1,0){.067156}}
\multiput(55.26,41.4)(.067156,-.03296){13}{\line(1,0){.067156}}
\multiput(57.01,40.54)(.067156,-.03296){13}{\line(1,0){.067156}}
\multiput(58.76,39.69)(.067156,-.03296){13}{\line(1,0){.067156}}
\multiput(60.5,38.83)(.067156,-.03296){13}{\line(1,0){.067156}}
\multiput(62.25,37.97)(.067156,-.03296){13}{\line(1,0){.067156}}
\multiput(63.99,37.12)(.067156,-.03296){13}{\line(1,0){.067156}}
\multiput(65.74,36.26)(.067156,-.03296){13}{\line(1,0){.067156}}
\multiput(67.49,35.4)(.067156,-.03296){13}{\line(1,0){.067156}}
\multiput(69.23,34.55)(.067156,-.03296){13}{\line(1,0){.067156}}
\multiput(70.98,33.69)(.067156,-.03296){13}{\line(1,0){.067156}}
\multiput(72.72,32.83)(.067156,-.03296){13}{\line(1,0){.067156}}
\multiput(74.47,31.98)(.067156,-.03296){13}{\line(1,0){.067156}}
\multiput(76.22,31.12)(.067156,-.03296){13}{\line(1,0){.067156}}
%\end
\put(48.37,44.87){\circle*{1}}
\put(77.25,30.75){\circle*{1}}
\put(91.12,61.62){\circle*{1}}
\put(73,69){\circle*{1}}
\put(93.87,61.87){$U_+$}
\put(72,70.75){$U_m$}
{\tiny
\put(78.37,45.87){$\overrightarrow{R}_2(U_+)$}
\put(63.77,53.22){$\overrightarrow{R}_2(U_m)$}
\put(82,38){$S_0(U_1,U_0)$}}
\qbezier(31.02,22.8)(78.84,32.75)(84.94,78.93)
\put(90.83,71.65){${\rm I}$}
\put(59.31,72.24){${\rm I\!I}$}
\put(42.22,71.5){${\rm I\!I\!I}$}
\qbezier(91.12,61.69)(86.14,42.66)(77.3,30.77)
%\dashline{1}(48.46,44.89)(43.85,54.41)
\multiput(48.39,44.82)(-.032225,.066529){13}{\line(0,1){.066529}}
\multiput(47.55,46.55)(-.032225,.066529){13}{\line(0,1){.066529}}
\multiput(46.71,48.28)(-.032225,.066529){13}{\line(0,1){.066529}}
\multiput(45.88,50.01)(-.032225,.066529){13}{\line(0,1){.066529}}
\multiput(45.04,51.74)(-.032225,.066529){13}{\line(0,1){.066529}}
\multiput(44.2,53.47)(-.032225,.066529){13}{\line(0,1){.066529}}
%\end
\qbezier(43.85,54.26)(38.65,49.95)(31.07,48.61)
\put(48,39){$U_c$}
\put(81.01,29.58){$U_{c*}$}
\put(43.85,54.11){\circle*{1}}
\put(44,56.78){$U_c^{*}$}
\put(27,48.31){$\tilde{h}$}
\put(55.3,18){$\bar{u}_+$}
\put(26.5,36.27){$\bar{h}_m$}
\put(56.93,79.38){$\Gamma_+$}
\put(87.41,79.38){$\Gamma_1$}
\end{picture}

%TeXCAD Picture [traffic flow 3'.pic]. Options:
%\grade{\on}
%\emlines{\off}
%\epic{\off}
%\beziermacro{\on}
%\reduce{\on}
%\snapping{\off}
%\quality{8.00}
%\graddiff{0.01}
%\snapasp{1}
%\zoom{11.3136}
\unitlength 0.7mm % = 2.85pt
\linethickness{0.4pt}
\ifx\plotpoint\undefined\newsavebox{\plotpoint}\fi % GNUPLOT compatibility
\begin{picture}(108.75,83.75)(32,10)
\footnotesize
\put(26.5,22.75){\vector(1,0){80}}
\put(31,19.75){\vector(0,1){64}}
\put(108.75,22.5){$u$}
\put(36.5,10.25){Fig. 3.5. $U_+\in {\rm I\!I}$, $U_m\in {\rm I\!I}$.}
\qbezier(31.01,22.7)(66.27,36.16)(71.9,82.83)
\qbezier(31.01,22.7)(93.55,32.53)(98.49,82.93)
\put(27.5,82){$h$}
\put(66,79){$\Gamma_+$}
\put(91,79){$\Gamma_1$}
%\dashline{1}(77.69,68.5)(84.32,56.39)
\multiput(77.62,68.43)(.031567,-.057663){14}{\line(0,-1){.057663}}
\multiput(78.51,66.82)(.031567,-.057663){14}{\line(0,-1){.057663}}
\multiput(79.39,65.2)(.031567,-.057663){14}{\line(0,-1){.057663}}
\multiput(80.28,63.59)(.031567,-.057663){14}{\line(0,-1){.057663}}
\multiput(81.16,61.97)(.031567,-.057663){14}{\line(0,-1){.057663}}
\multiput(82.04,60.36)(.031567,-.057663){14}{\line(0,-1){.057663}}
\multiput(82.93,58.74)(.031567,-.057663){14}{\line(0,-1){.057663}}
\multiput(83.81,57.13)(.031567,-.057663){14}{\line(0,-1){.057663}}
%\end
\put(85.56,56.83){$U_+$}
\put(79.2,68.77){$U_m$}
\qbezier(77.78,68.41)(67.93,53.74)(31.02,46.14)
\qbezier(84.23,56.48)(62.4,41.98)(31.02,34.56)
{\tiny
\put(78.93,46.85){$S_0(U_1,U_0)$}
\put(49.41,51.71){$\overrightarrow{R}_2(U_m)$}
\put(60,44.11){$\overrightarrow{R}_2(U_+)$}}
\qbezier(84.23,56.48)(77.87,46.45)(72.92,43.66)
%\dashline{1}(72.92,43.66)(65.94,58.07)
\multiput(72.85,43.59)(-.031596,.065192){13}{\line(0,1){.065192}}
\multiput(72.03,45.29)(-.031596,.065192){13}{\line(0,1){.065192}}
\multiput(71.21,46.98)(-.031596,.065192){13}{\line(0,1){.065192}}
\multiput(70.39,48.68)(-.031596,.065192){13}{\line(0,1){.065192}}
\multiput(69.56,50.37)(-.031596,.065192){13}{\line(0,1){.065192}}
\multiput(68.74,52.07)(-.031596,.065192){13}{\line(0,1){.065192}}
\multiput(67.92,53.76)(-.031596,.065192){13}{\line(0,1){.065192}}
\multiput(67.1,55.46)(-.031596,.065192){13}{\line(0,1){.065192}}
\multiput(66.28,57.15)(-.031596,.065192){13}{\line(0,1){.065192}}
%\end
\put(65.94,57.98){\circle*{1}}
\put(72.92,43.66){\circle*{1}}
\put(100,64.22){${\rm I}$}
\put(90,64.64){${\rm I\!I}$}
\put(42.22,71.5){${\rm I\!I\!I}$}
%\dashline{1}(65.94,57.98)(57.81,67.79)
\multiput(65.87,57.91)(-.032269,.038933){18}{\line(0,1){.038933}}
\multiput(64.71,59.31)(-.032269,.038933){18}{\line(0,1){.038933}}
\multiput(63.54,60.72)(-.032269,.038933){18}{\line(0,1){.038933}}
\multiput(62.38,62.12)(-.032269,.038933){18}{\line(0,1){.038933}}
\multiput(61.22,63.52)(-.032269,.038933){18}{\line(0,1){.038933}}
\multiput(60.06,64.92)(-.032269,.038933){18}{\line(0,1){.038933}}
\multiput(58.9,66.32)(-.032269,.038933){18}{\line(0,1){.038933}}
%\end
\qbezier(58.51,67)(49.1,59.66)(31.02,55.15)
\put(58.43,67){\circle*{1}}
\put(68.24,57.28){$U_c$}
\put(71,39){$U_{c*}$}
\put(60,68.15){$U_c^{*}$}
\put(26,55.15){$\tilde{h}$}
\put(26.61,46.05){$\bar{h}_m$}
\put(26,34.29){$\bar{h}_+$}
\end{picture}
%TeXCAD Picture [traffic flow 5.pic]. Options:
%\grade{\on}
%\emlines{\off}
%\epic{\off}
%\beziermacro{\on}
%\reduce{\on}
%\snapping{\off}
%\quality{8.00}
%\graddiff{0.01}
%\snapasp{1}
%\zoom{9.5137}
\unitlength 0.7mm % = 2.85pt
\linethickness{0.4pt}
\ifx\plotpoint\undefined\newsavebox{\plotpoint}\fi % GNUPLOT compatibility
\begin{picture}(50.75,50.75)(44,10)
\footnotesize
\put(26.5,22.75){\vector(1,0){80}}
\put(31,19.75){\vector(0,1){64}}
\put(108.75,22.5){$u$}
\put(27,82){$h$}
\put(36.5,10.25){Fig. 3.6. $U_+\in {\rm I\!I\!I}$, $U_m\in {\rm I\!I\!I}$.}
\qbezier(31.01,22.7)(84.93,28.38)(93.02,79.04)
\qbezier(31.01,22.7)(70.74,37.52)(73.89,78.83)
\put(67,76){$\Gamma_+$}
\put(87,76){$\Gamma_1$}
%\dashline{1}(58.23,67.38)(62.33,53.4)
\multiput(58.16,67.31)(.03037,-.10355){9}{\line(0,-1){.10355}}
\multiput(58.71,65.44)(.03037,-.10355){9}{\line(0,-1){.10355}}
\multiput(59.25,63.58)(.03037,-.10355){9}{\line(0,-1){.10355}}
\multiput(59.8,61.71)(.03037,-.10355){9}{\line(0,-1){.10355}}
\multiput(60.35,59.85)(.03037,-.10355){9}{\line(0,-1){.10355}}
\multiput(60.89,57.99)(.03037,-.10355){9}{\line(0,-1){.10355}}
\multiput(61.44,56.12)(.03037,-.10355){9}{\line(0,-1){.10355}}
\multiput(61.99,54.26)(.03037,-.10355){9}{\line(0,-1){.10355}}
%\end
\qbezier(58.23,67.27)(48.93,59.55)(31.01,55.81)
\qbezier(62.23,53.61)(51.56,45.67)(31.01,40.05)
{\tiny
\put(41.73,61){$\overrightarrow{R}_2(U_+)$}
\put(43.73,44.46){$\overrightarrow{R}_2(U_m)$}
\put(42,53){$S_0(U_1,U_0)$}}
\put(60.12,67.48){$U_+$}
\put(63.59,53.5){$U_m$}
\put(93.86,64.22){${\rm I}$}
\put(83.14,64.64){${\rm I\!I}$}
\put(42.22,71.5){${\rm I\!I\!I}$}
\qbezier(58.23,67.27)(51.03,56.44)(31.01,48.98)
\put(26.5,55.5){$\bar{h}_+$}
\put(26.5,39.84){$\bar{h}_m$}
\put(26.5,48.77){$\tilde{h}$}
\end{picture}\\
\
\noindent
\begin{proof}
From \eqref{3.2}, one has $U_0\in \overrightarrow{R}_2(U_m,U)$, which is
\begin{equation}\label{3.4}
 u_0=u_m+2\sqrt{g}(\sqrt{h_0}-\sqrt{h_m}),\quad 0<h_0<h_m.
\end{equation}
$\overrightarrow{R}_2(U_+,U)$ is parameterized by $h_1$ as
\begin{equation}\label{3.5}
 u=u(h_1)=u_++2\sqrt{g}(\sqrt{h_1}-\sqrt{h_+}),\quad 0\leq h_1 \leq h_+.
\end{equation}
Besides, $U_0$ and $U_1$ are connected by a stationary wave
\begin{equation}\label{3.6}
\left\{
\begin{array}{ll}
h_0u_0=h_1u_1,\\[5pt]
\displaystyle \frac{u_0^2}2+g(h_0+a_0)=\frac{u_1^2}2+g(h_1+a_1).
\end{array}
\right.
\end{equation}
Differentiating \eqref{3.6} one has
\begin{equation}\label{3.7}
\left\{
\begin{array}{ll}
u_0{\rm d}h_0+h_0{\rm d}u_0=u_1{\rm d}h_1+h_1{\rm d}u_1,\\[5pt]
u_0{\rm d}u_0+g{\rm d}h_0=u_1{\rm d}u_1+g{\rm d}h_1.
\end{array}
\right.
\end{equation}
We note \eqref{3.4} yields 
\begin{equation}\label{3.8}
\frac{{\rm d}u_0}{{\rm d}h_0}=\sqrt{\frac{g}{h_0}}.
\end{equation}
Substituting \eqref{3.8} into \eqref{3.7}, we get the derivative of $u_1(h_1;U_0)$ as
\begin{equation}\label{3.9}
\displaystyle \frac{{\rm d}u_1}{{\rm d}h_1}=\sqrt{\frac{g}{h_0}} \cdot \frac{u_1-\sqrt{g h_0}}{u_1-\sqrt{gh_1/h_0}}=\frac{\sqrt{g}u_1-g\sqrt{h_0}}{\sqrt{h_0}u_1-\sqrt{g}h_1}.
\end{equation}
As $S_0(U_1,U_0)$ can be written as a parameterized function, denoted by $u_1(h_1;U_0)$,
we next prove that $u_1(h_1;U_0)$ is on the right of $\overrightarrow{R}_2(U_+,U)$ (which is parameterized as $u(h_1)$) in the same domain $0\leq h_1\leq h_+$. Defining
\begin{equation}\label{3.10}
f(h_1)=u_1(h_1)-u(h_1),\quad 0\leq h_1\leq h_+,
\end{equation}
one has $f(h_+)=0$. Moreover, combining \eqref{3.5} and \eqref{3.9}, one gets
\begin{equation}\label{3.11}
\frac{{\rm d}f}{{\rm d}h_1}=\frac{{\rm d}u_1}{{\rm d}h_1}-\frac{{\rm d}u}{{\rm d}h_1}=\sqrt{\frac{g}{h_1}}\cdot \frac{(u_1+\sqrt{gh_1})(\sqrt{h_1}-\sqrt{h_0})}{u_1\sqrt{h_0}-\sqrt{g}h_1}.
\end{equation}

Similarly, denote $q(u_1)=h_1(u_1)-h(u_1),\quad 0\leq u_1\leq u_+$, where $h(u_1)$ is solved in \eqref{3.5}
\begin{equation}\label{3.12}
\displaystyle h=h(u_1)=\frac{(u_1-u_++2\sqrt{gh_+})^2}{4g}.
\end{equation}
Besides, the monotone property of $h_1(u_1)$ is obtained from \eqref{3.7},
\begin{equation}\label{3.13}
\displaystyle \frac{{\rm d}h_1}{{\rm d}u_1}=\frac{u_1\sqrt{h_0}-\sqrt{g}h_1}{u_1\sqrt{g}-g\sqrt{h_0}}.
\end{equation}
Finally, one gets
\begin{equation}\label{3.14}
\displaystyle \frac{{\rm d}q}{{\rm d}u_1}=\frac{{\rm d}h_1}{{\rm d}u_1}-\frac{{\rm d}h}{{\rm d}u_1}=\frac{(u_1+\sqrt{gh_1})(\sqrt{h_0}-\sqrt{h_1})}{\sqrt{g}(u_1-\sqrt{gh_0})}.
\end{equation}

Following the assumption $a_0>a_1$, it leads to
\begin{equation}\label{3.15}
u_1\sqrt{h_0}-h_1\sqrt{g}\left\{
\begin{array}{ll}
>0\quad U_1\in {\rm I}\\
<0\quad U_1\in {\rm I\!I \cup I\!I\!I},
\end{array}\right.
\quad {\rm and}\quad 
u_1-\sqrt{gh_0}\left\{
\begin{array}{l}
>0\quad U_1\in {\rm I}\\
<0\quad U_1\in {\rm I\!I \cup I\!I\!I}.
\end{array}\right.
\end{equation}
We have consequently that
\begin{equation}\label{3.16}
\frac{{\rm d}f}{{\rm d}h_1}\left\{
\begin{array}{ll}
<0\quad U_1\in {\rm I\!I \cup I\!I\!I}\\
<0\quad U_1\in {\rm I},
\end{array}\right.
\quad {\rm and}\quad 
\frac{{\rm d}q}{{\rm d}u_1}\left\{
\begin{array}{l}
>0\quad U_1\in {\rm I\!I \cup I\!I\!I}\\
>0\quad U_1\in {\rm I}.
\end{array}\right.
\end{equation}
Therefore, we make the following conclusions.

\noindent
1. $U_+\in {\rm I}$, $U_m\in {\rm I}$. $S_0(U_1,U_0)$ intersects with $u-$axis at $(\tilde{u},0)$, since $f'(h_1)<0$ and $f(h_+)=0$, we have $f(h_1)<0$ when $0\leq h_1\leq h_+$. Thus $S_0(U_1,U_0)$ is on the right of $\overrightarrow{R}_2(U_+,U)$. See Fig. 3.3.

\noindent
2. $U_+\in {\rm I}$, $U_m\in {\rm I\!I}$. Denote $U_c=\overrightarrow{R}_2(U_m,U)\cap \Gamma_+$. A similar discussion shows that $S_0(U_1,U_0)$ is on the right of $\overrightarrow{R}_2(U_+,U)$ in the supercritical area. In the subcritical area, $S_0$ is on the left of $\overrightarrow{R}_2(U_+,U)$. See Fig. 3.4.

\noindent
3. $U_+\in {\rm I\!I}$, $U_m\in {\rm I\!I}$. This case is similar to case 2. See Fig. 3.5.

\noindent
4. $U_+\in {\rm I\!I\!I}$, $U_m\in {\rm I\!I\!I}$. $g(u_1)<0$ when $0\leq u_1\leq u_+$. $S_0(U_1,U_0)$ is below $\overrightarrow{R}_2(U_+,U)$. See Fig. 3.6.
\end{proof}
When $U_+\in {\rm I}$, $U_m\in {\rm I}$, denote $\overrightarrow{R}_2(U_+)\cap \{h=0\}=(\bar{u}_+,0)$, $\overrightarrow{R}_2(U_m)\cap \{h=0\}=(\bar{u}_m,0)$, $S_0(U_1,U_0)\cap \{h=0\}=(\tilde{u},0)$, see Fig. 3.3. Next it can be shown that $\bar{u}_+$ is indeed on the left side of $\tilde{u}$ as indicated in lemma \ref{3.1}. In fact, $(\tilde{u},0)$ and $(\bar{u}_m,0)$ are connected by a stationary wave $S_0$. From \eqref{2.14}, one has
\begin{equation}\label{3.17}
\displaystyle \frac{\tilde{u}^2}2+ga_1=\frac{\bar{u}_m^2}2+ga_0.
\end{equation}
Besides, since $(\bar{u}_+,0)\in \overrightarrow{R}_2(U_+,U)$ and $(\bar{u}_m,0)\in \overrightarrow{R}_2(U_m,U)$, one gets
 \begin{equation}\label{3.18}
\displaystyle \bar{u}_+=u_+-2\sqrt{g h_+}, \quad \displaystyle \bar{u}_m=u_m-2\sqrt{g h_m}.
\end{equation}
A direct calculation shows that
 \begin{equation}\label{3.19}
\displaystyle \tilde{u}^2-\bar{u}_+^2=2g(h_m-h_+)+4\sqrt{g_+u_+(u_+-u_m)}>0,
\end{equation}
which follows from the fact $\tilde{u}>\bar{u}_+$.

Similarly, when $U_+\in {\rm I\!I\!I}$, $U_m\in {\rm I\!I\!I}$, denote $\overrightarrow{R}_2(U_+)\cap \{u=0\}=(\bar{h}_+,0)$, $\overrightarrow{R}_2(U_m)\cap \{u=0\}=(\bar{h}_m,0)$, $S_0(U_1,U_0)\cap \{u=0\}=(\tilde{h},0)$.
See Fig. 3.6, one may see that $\tilde{h}<\bar{h}_+$. In fact, from the above discussion, we have
 \begin{equation}\label{3.20}
\displaystyle \tilde{h}+a_1=\bar{h}_m+a_0,
\end{equation}
and
\begin{equation}\label{3.21}
\left\{
\begin{array}{ll}
\displaystyle u_+-2\sqrt{gh_+}=-2\sqrt{g\bar{h}_+},\\[5pt]
\displaystyle u_m-2\sqrt{gh_m}=-2\sqrt{g\bar{h}_m}.
\end{array}
\right.
\end{equation}
It follows that
 \begin{equation}\label{3.22}
 \begin{array}{ll}
  \displaystyle \tilde{h}-\bar{h}_+ & =\displaystyle \frac{h_m-h_+}2+\frac{a_0-a_1}{2}+\frac{\sqrt{h_+u_+}(\sqrt{u_+}-\sqrt{u_m})}{\sqrt{g}}\\[7pt]
 & \displaystyle =\frac{u_+^2-u_m^2}{2}+\frac{\sqrt{h_+u_+}(\sqrt{u_+}-\sqrt{u_m})}{\sqrt{g}}<0.
\end{array}
\end{equation}

Now we discuss the interaction results case by case.

\noindent
{\bf Case 1.} $U_+\in {\rm I}$, $U_m\in {\rm I}$. See Fig. 3.7.
In this case, from the result 1) in lemma \ref{lem3.2}, we draw the curves $\overrightarrow{R}_2(U_m,U), \overrightarrow{R}_2(U_+,U)$ and $S_0(U_1,U_0)$, which intersect $u-$axis at $\bar{u}_m, \bar{u}_+$ and $\tilde{u}$ respectively. Denote $U_{-*}\in S_0(U_1,U_0)$, which is obtained by $U_-\in \overrightarrow{R}_2(U_m,U)$. Denote $U_2=\overleftarrow{S}_1(U,U_{-*})\cap \overrightarrow{R}_2(U_+,U)$.

\begin{figure}[htbp]
\subfigure{
\begin{minipage}[t]{0.45\textwidth}
\centering
\includegraphics[width=\textwidth]{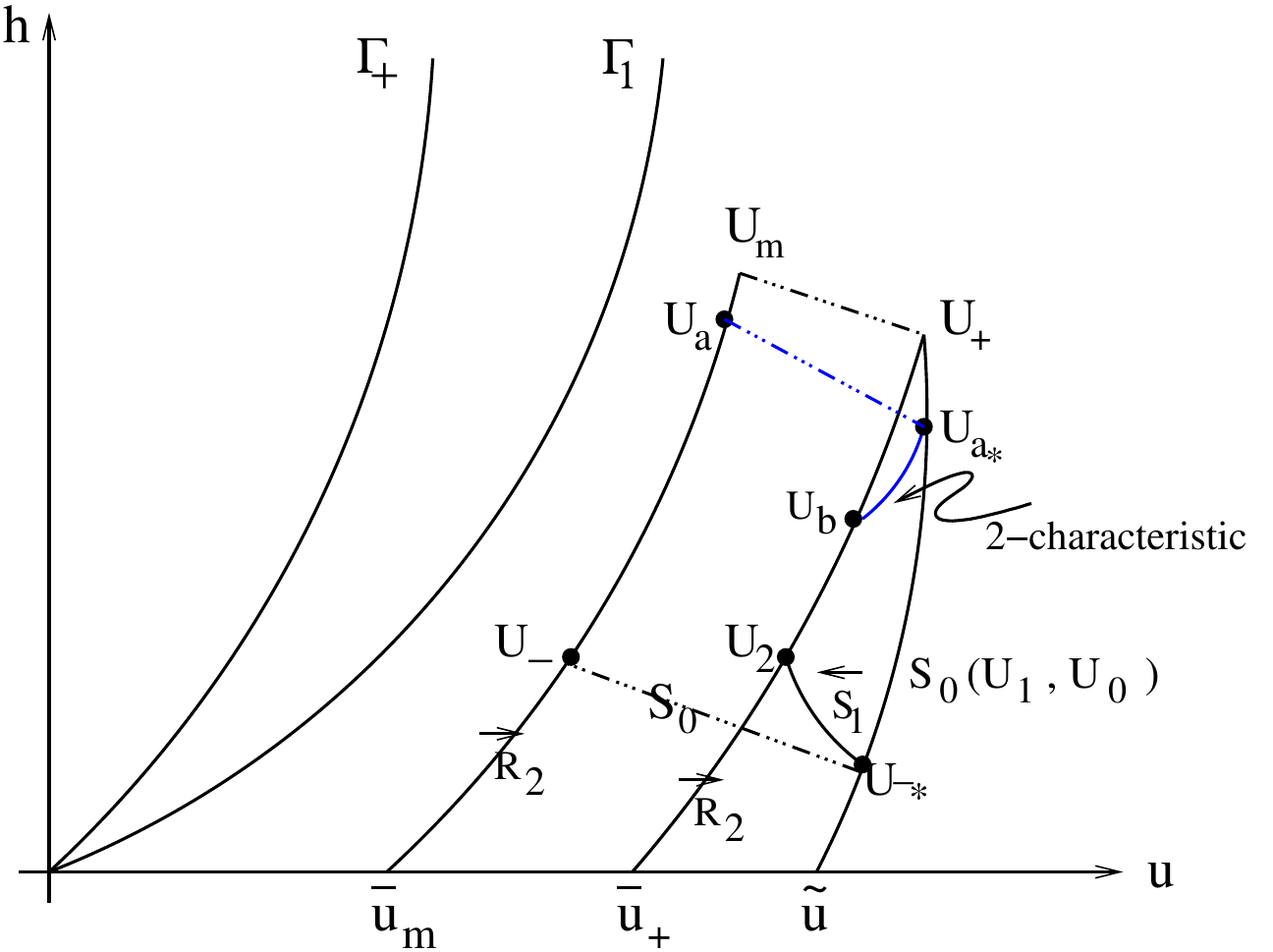}
\end{minipage}
}
\subfigure{
\begin{minipage}[t]{0.45\textwidth}
\centering
\includegraphics[width=\textwidth]{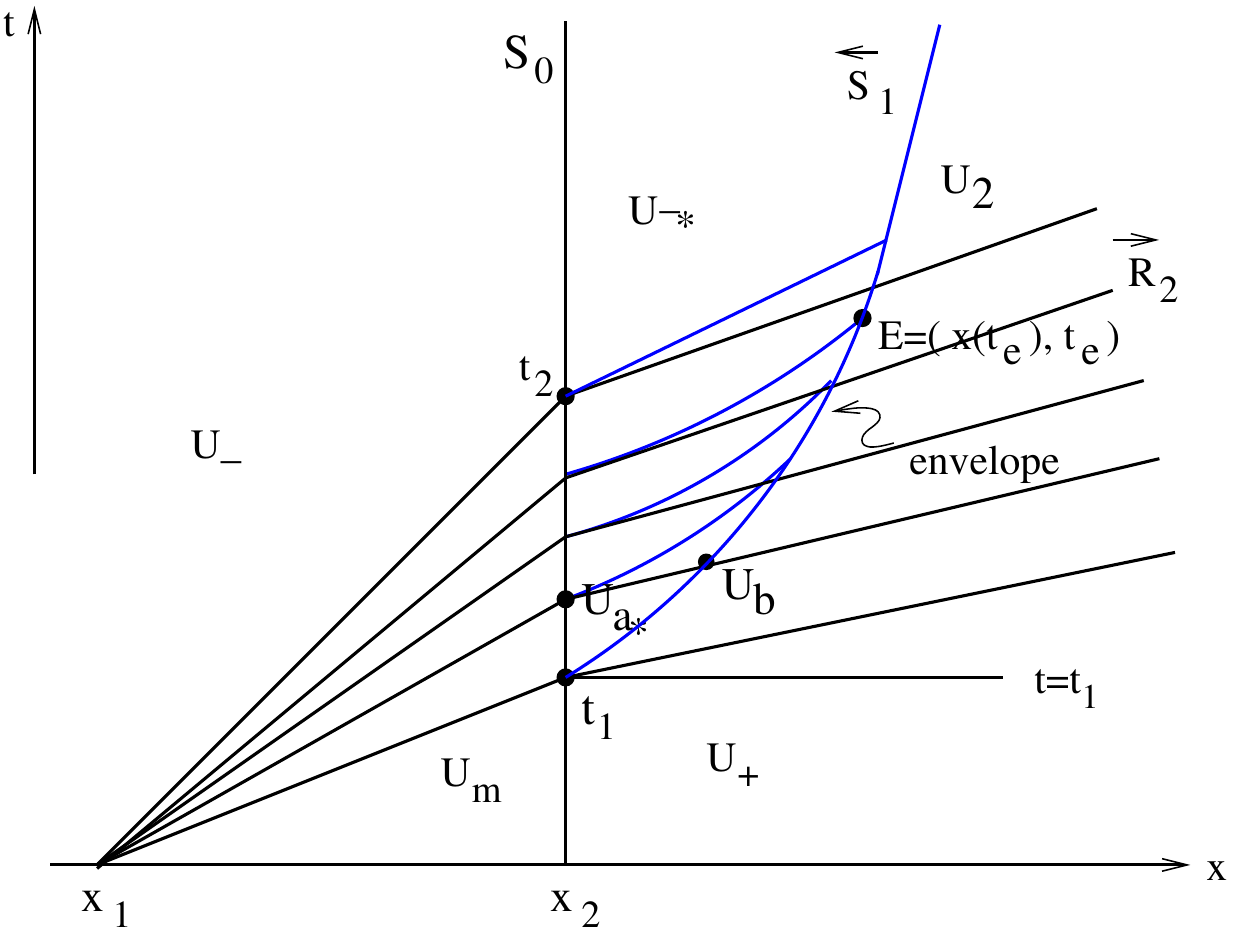}
\end{minipage}
}
\caption*{Fig. 3.7. Case 1. $U_+\in {\rm I}$, $U_m\in {\rm I}$.}
\end{figure}

To determine that the 1-characteristic curves are compressive during the interaction process, we denote $U_a\in \overrightarrow{R}_2(U_m,U)$, $U_a$ jumps to $U_{a*}$ by a stationary wave $S_0$. From $U_{a*}$, we draw a 2-characteristic curve which intersects with $\overrightarrow{R}_2(U_+,U)$ at $U_b$, then one has $u_{a*}>u_b$, see Fig. 3.7 (left). Since the velocity $u$ increases across the wave $\overrightarrow{R}_2$ from back side, we conclude that the 1-characteristics are compressible in the interaction domain. See Fig. 3.7 (right).

The compressible 1-characteristics will form an envelope, which will develop to a shock wave $\overleftarrow{S}_1$ as time goes on. Denote the starting point of the shock wave  
$\overleftarrow{S}_1$ as $E=(x(t_e),t_e)$. $E$ might lie either in or out of the interaction area. We discuss the results as follows.

\noindent%
1). When $E$ is out of the interaction area, we know that the shock wave propagates on the right with a constant speed.

\noindent%
2). When $E$ is in the interaction area, we solve a free boundary value problem of  the 1-shock wave $\overleftarrow{S}_1$ at $(x(t_e),t_e)$
\begin{equation}\label{3.23}
\left\{
\begin{array}{lll}
\overleftarrow{S}_1(u_r,u_l):
\left\{
\begin{array}{lll}
\displaystyle \frac{{\rm d}x}{{\rm d}t}=u_l-\sqrt{\frac{g}2(h_r+h_l)\frac{h_r}{h_l}},\\[8pt]
\displaystyle u_r=u_l-\sqrt{\frac{g}2(h_r-h_l)^2(\frac1{h_r}+\frac1{h_l})},  h_r>h_l. \\[8pt]
x(t_e)=x_e,
\end{array}\right.~ (\rm{the~RH~condition~of~\overleftarrow{S}_1})\\[7pt]
\displaystyle  \overrightarrow{R}_2(U,U_2): {\rm d}u_r=\frac{g}{h_r}, h_2<h<h_+.
~\qquad   \rm{(U_r~\in \overrightarrow{R}_2(U,U_2)~on~the~right~of~the~shock)}\\[10pt]
(u,h)=\left\{
\begin{array}{ll}
(u_+,h_+),\quad 0<t<t_1,~x=x_2.\\
(u_1,h_1),~U_1\in S_0(U_1,U_0) ~{\rm with}~U_0\in \overleftarrow{R}_1(U,U_m),~t_1<t<t_2,~x=x_2,\\[5pt]
\displaystyle (u_{-*},h_{-*}),\quad \frac{x-x_2}{t-t_2}=u_{-*}+c_{-*},~x>x_2
\end{array}
\right.
\end{array}
\right.
\end{equation}
in the rectangular domain $\{(x,t)\big|x\geq x_e, t\geq t_e\}$, where $x=x(t)$ is the shock wave supplemented by the Lax entropy condition: $\displaystyle 0<u_r-2\sqrt{gh_r}<  \frac{{\rm d}x}{{\rm d}t}<u_l-2\sqrt{gh_l}$. Thus the shock wave $\overleftarrow{S}_1$ propagates with a positive speed in the domain. From \cite{ChangHsiao3}, it is known the speed of $\overleftarrow{S}_1$ will decrease during the process of penetrating $\overrightarrow{R}_2$, see Fig. 3.7 (left). Finally, the large time behavior of the solution is
\begin{equation}
U_-\oplus S_0(U_{-*},U_-)\oplus\overleftarrow{S}_1(U_2,U_{-*})\oplus \overrightarrow{R}_2(U_+,U_2)\oplus U_+,
\end{equation}
where $``\oplus"$ means ``follows".

\noindent
{\bf Case 2.} $U_+\in {\rm I}$, $U_m\in {\rm I\!I}$. See Fig. 3.8.

\noindent
{\bf Subcase 2.1.} As $U_-\notin {\rm I\!I\!I}$, which is on $\overrightarrow{R}_2(U_m,U)$ and in the supercritical area, it is similar to case 1. We omit it. See Fig. 3.7.

\noindent
{\bf Subcase 2.2.} As $U_-\in {\rm I\!I\!I}$, which is in the subcritical area. $S_0(U_1,U_0)$ touches $\Gamma_{*}$ at $U_{c*}$ and coincides with $\widehat{\widetilde{U}U_{c*}}$ on $\Gamma_{*}$ in the supercritical area. The interaction process includes two parts. In the first part from $t_1$ to $t_2$, $U_0\in \overrightarrow{R}_2(U_m,U)$ ($U_0$ is between $U_m$ and $U_c$), $U_1\in S_0(U_1,U_0)$ ($U_1$ is between $U_+$ and $U_{c*}$). This process is the same as the above case. In the second part from $t_2$ to $t_3$, $U_0\in \Gamma_+\cap \overleftarrow{R}_1(U,\overline{U}_m)$ ($\overline{U}_m\in \overrightarrow{R}_2(U_m,U)$ is between $U_c$ and $U_-$), $U_1\in S_0(U_1,U_0)$ ($U_1$ is between $\widetilde{U}=(\tilde{u},0)$ and $U_{c*}$).

On the right of $x=x_2$, one can solve a free boundary value problem as \eqref{3.23} at the starting of the shock $\overleftarrow{S}_1$ at $(x(t_e),t_e)$. While on the left of $x=x_2$, we solve an initial boundary value problem of \eqref{2.11} with
\begin{equation}\label{3.24}
(u,h)=\left\{
\begin{array}{l}
\displaystyle \big(\xi-\sqrt{gh}, \big(\frac{u_--2\sqrt{gh_-}-\xi}{3\sqrt{g}}\big)^2\big),~~~(x,t)\in C_-,~u_-+\sqrt{gh_-}<\xi<u_m+\sqrt{gh_m}, \\[8pt]
\displaystyle u=\sqrt{gh}, x=x_2,~t_2<t<t_3,\\[7pt]
\displaystyle (u_2,h_2),~x=x_2,t>t_3,
\end{array}
\right.
\end{equation}

\begin{figure}[htbp]
\subfigure{
\begin{minipage}[t]{0.45\textwidth}
\centering
\includegraphics[width=\textwidth]{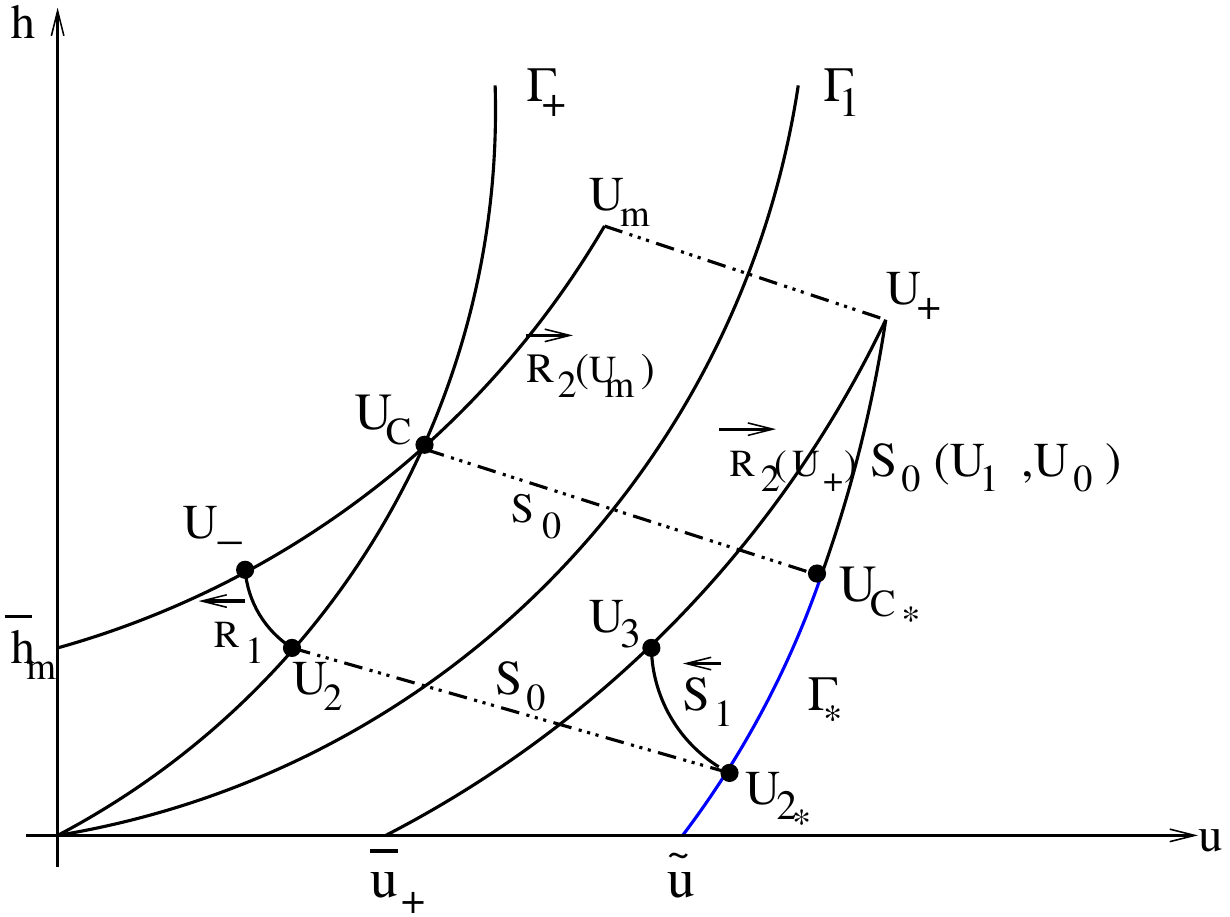}
\end{minipage}
}
\subfigure{
\begin{minipage}[t]{0.45\textwidth}
\centering
\includegraphics[width=\textwidth]{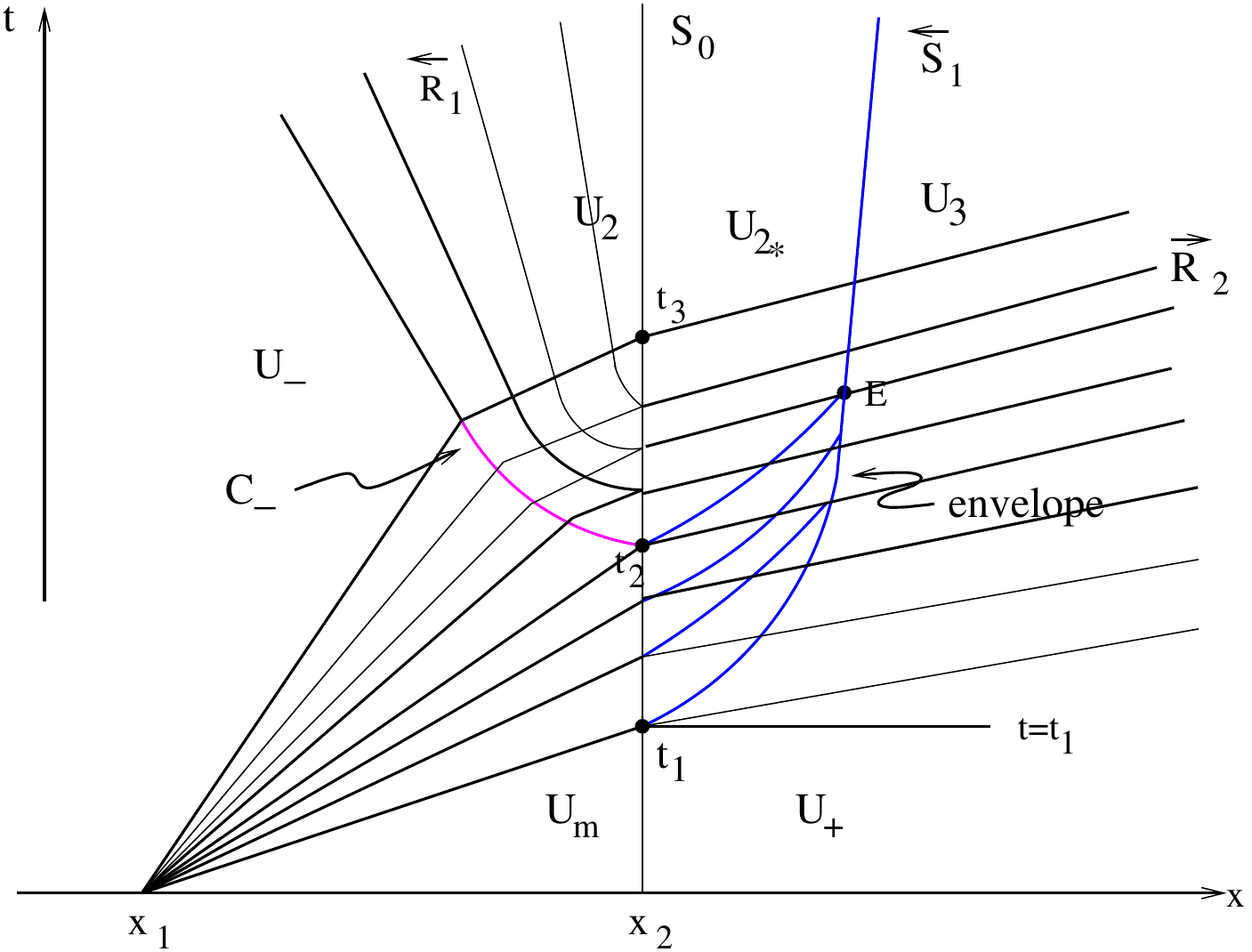}
\end{minipage}
}
\caption*{Fig. 3.8. Case 2. $U_+\in {\rm I}$, $U_m\in {\rm I\!I}$.}
\end{figure}

\noindent%
where $\xi$ is the given slope of the characteristic line of $\overleftarrow{R}_1(U_m,U_-)$, $C_-$ denotes the penetrating backward characteristic, see Fig. 3.8 (right). The existence and uniqueness of the two problems can be obtained by the classical theory from Li and Yu (\cite{LiTT}), see also in Wang and Wu (\cite{Wang}). Furthermore, the large time behavior of the solution from the theory of $p-$system \eqref{2.11} is
\begin{equation}
U_-\oplus \overleftarrow{R}_1(U_2,U_-)\oplus S_0(U_{2*},U_2)\oplus \overleftarrow{S}_1(U_3,U_{2*})\oplus \overrightarrow{R}_2(U_+,U_3)\oplus U_+.
\end{equation}

\begin{figure}[htbp]
\centering
\includegraphics[width=0.5\textwidth]{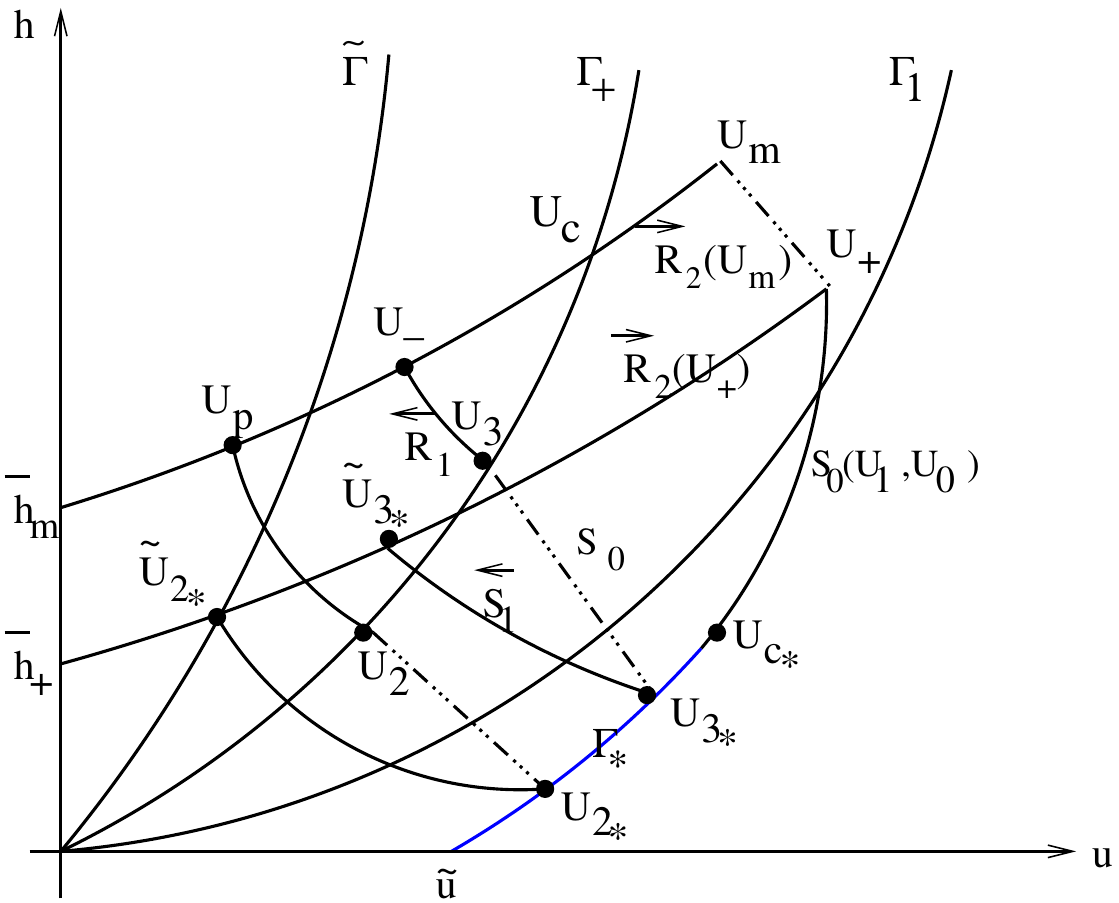}
\caption*{Fig. 3.9. Case 3. $U_+\in {\rm I\!I}, U_m\in {\rm I\!I}$.}
\end{figure}

\noindent
{\bf Case 3.} $U_+\in {\rm I\!I}, U_m\in {\rm I\!I}$. See Fig. 3.9.\\
From lemma \ref{lem3.2}, $S_0(U_1,U_0)$ is below $\overrightarrow{R}_2(U_+,U)$. It touches $\Gamma_{*}$ at $U_{c*}$ and coincides with $\widehat{\widetilde{U}U_{c*}}$ on $\Gamma_*$ in the supercritical area.

From any point $U\in \Gamma_*$, there exists a point $\widetilde{U}\in {\rm I\!I\!I}$,  such that the 1-shock wave speed vanishes, i.e., $\sigma(\widetilde{U},U)=0$.  Such states $\widetilde{U}$ form a curve in ${\rm I\!I\!I}$, denoted by $\widetilde{\Gamma}$, see Fig. 3.9.

Denote $\widetilde{U}_{2*}=\overrightarrow{R}_2(U_+,U)\cap \widetilde{\Gamma}$, which is obtained by $U_{2*}\in \Gamma_*$ with $\sigma(\widetilde{U}_{2*},U_{2*})=0$. $U_2\in \Gamma_+$ jumps to $U_{2*}$ by a stationary wave $S_0(U_{2*},U_2)$. Denote $U_p=\overleftarrow{R}_1(U_2,U)\cap \overrightarrow{R}_2(U_m,U)$.

By virtue of \eqref{2.22}, we know that $U_p$ is a critical point. When $U_-$ is above $U_p$ on $\overrightarrow{R}_2(U_m,U)$, $U_-$ first reaches to $U_3$ on $\Gamma_+$ by $\overleftarrow{R}$, $U_3$ jumps to $U_{3*}$ by a stationary wave $S_0$, then $U_{3*}$ reaches to $\widetilde{U}_{3*}$ by $\overleftarrow{S}_1(\widetilde{U}_{3*},U_{3*})$ with a positive speed, see Fig. 3.9.  As soon as $U_-$ is below $U_p$ on $\overrightarrow{R}_2(U_m,U)$, the 1-shock wave $\overleftarrow{S}_1$ has a negative speed.

\noindent
{\bf Subcase 3.1.} $U_-$ is between $U_m$ and $U_c$ on the curve $\overrightarrow{R}_2(U_m,U)$, it is similar to case 1, see Fig. 3.7. We omit here.

\noindent
{\bf Subcase 3.2.} $U_-$ is between $U_c$ and $U_p$ on the curve $\overrightarrow{R}_2(U_m,U)$. The interaction process includes two parts. The first part is from time $t_1$ to $t_2$, $U_0$ is between $U_m$ and $U_c$. A 1-shock wave with positive speed emits in this part. The second part is from $t_2$ to $t_3$, $U_0$ is between $U_c$ and $U_p$. A backward rarefaction wave propagates during this process. Both the two parts are the same with case 2, see Fig. 3.8.

We note the speed of the 1-shock wave $\overleftarrow{S}_1$ decreases in
$\overrightarrow{R}_2$ as discussed. At the critical case $U_-=U_p$, where $U_3=U_2$, $U_{3*}=U_{2*}$, the speed of the shock is equal to zero, and then becomes negative as $U_-$ is below $U_p$ on $\overrightarrow{R}_2(U_m,U)$. See the following case.

\begin{figure}[htbp]
\centering
\includegraphics[width=0.5\textwidth]{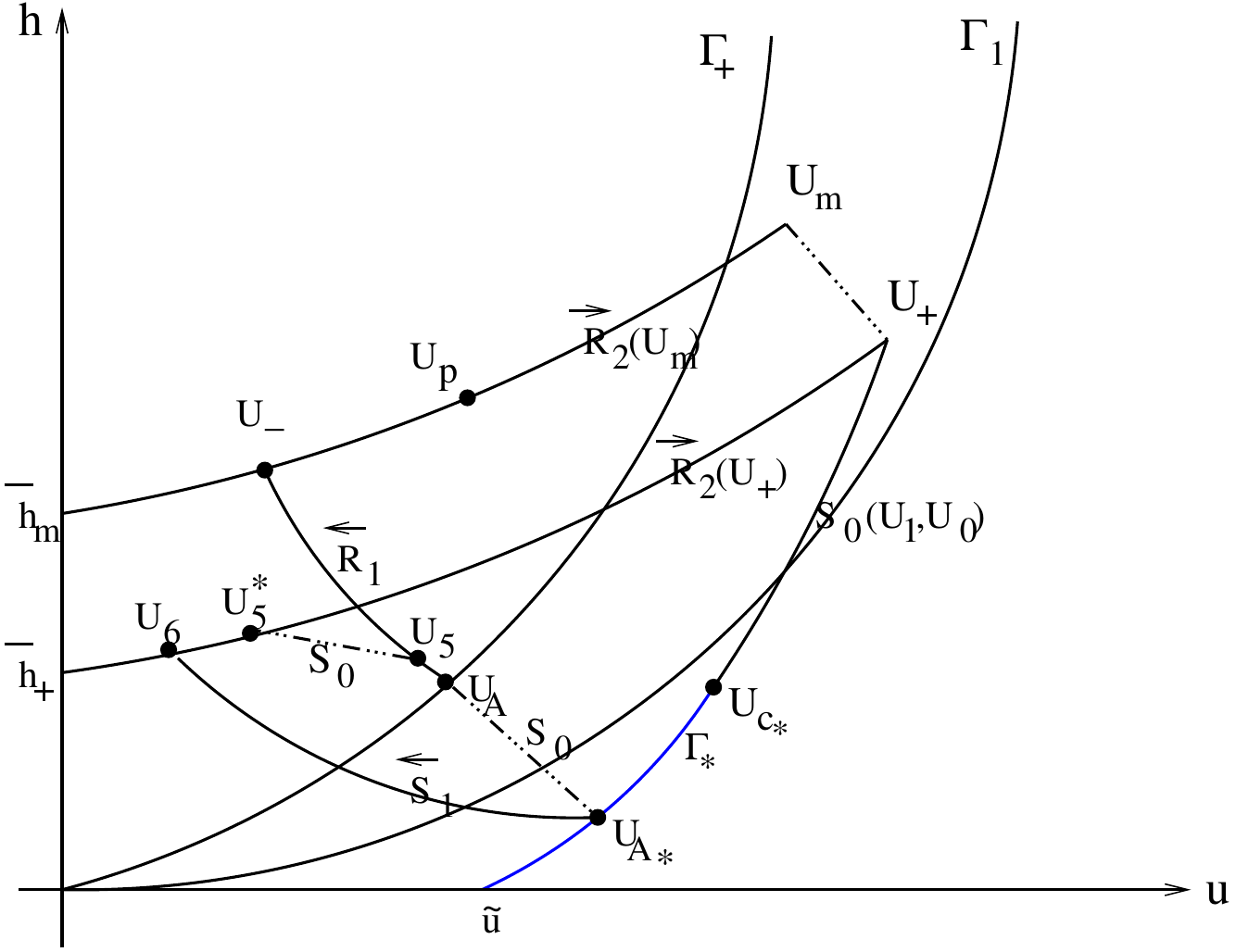}
\caption*{Fig. 3.10. Subcase 3.3. $U_-$ is between $U_p$ and $\overline{U}_m$.}
\end{figure}

\noindent
{\bf Subcase 3.3.} $U_-$ is between $U_p$ and $\overline{U}_m$ on the curve $\overrightarrow{R}_2(U_m,U)$. The interaction includes three parts. The first two parts from time $t_1$ to $t_3$ are the same with subcase 3.2. When $U_0$ goes down along $\Gamma_+$ passing $U_p$, a shock wave $\overleftarrow{S}_1$ transmits during the interaction. The speed of $\overleftarrow{S}_1$ decreases from time $t_1$ to $t_3$. At the critical case $U_0=U_2$, $U_1=U_{2*}$, the speed equals to zero and then becomes negative when $U_0$ is below $U_p$ on the curve $\Gamma_+$. It propagates on the left with a constant negative speed before touching $S_0$. See Fig. 3.11 (left).

The third part is from time $t_3$ to $t_4$, $\overleftarrow{S}_1$ penetrates $\overrightarrow{R}_2$ and touches $S_0$ during this period. When $\overleftarrow{S}_1$ touches $S_0$ at $t=t_4$, we solve a generalized Riemann problem in this case. The initial data on the left side of $S_0$ is on $\overleftarrow{R}_1(U,U_-), u\geq u_-$, the right side of $S_0$ is a constant state $U_6$. $U_6$ is connected with $U_{A*}$ by $\overleftarrow{S}_1$. See Fig. 3.10.

By solving the generalized Riemann problem, the solution includes a  backward rarefaction wave $\overleftarrow{R}_1$ and $\overleftarrow{S}_1$, which connects $U_-$ and $U_5$, $U_5$ jumps to $U_{5}^{*}$ by a stationary wave $S_0(U_5^{*},U_5)$, finally $U_5^{*}$ connects with $U_6$ by $\overrightarrow{S}_2(U_6,U_5^{*})$. See Fig. 3.11 (left).

\begin{figure}[htbp]
\subfigure{
\begin{minipage}[t]{0.47\textwidth}
\centering
\includegraphics[width=1.0\textwidth]{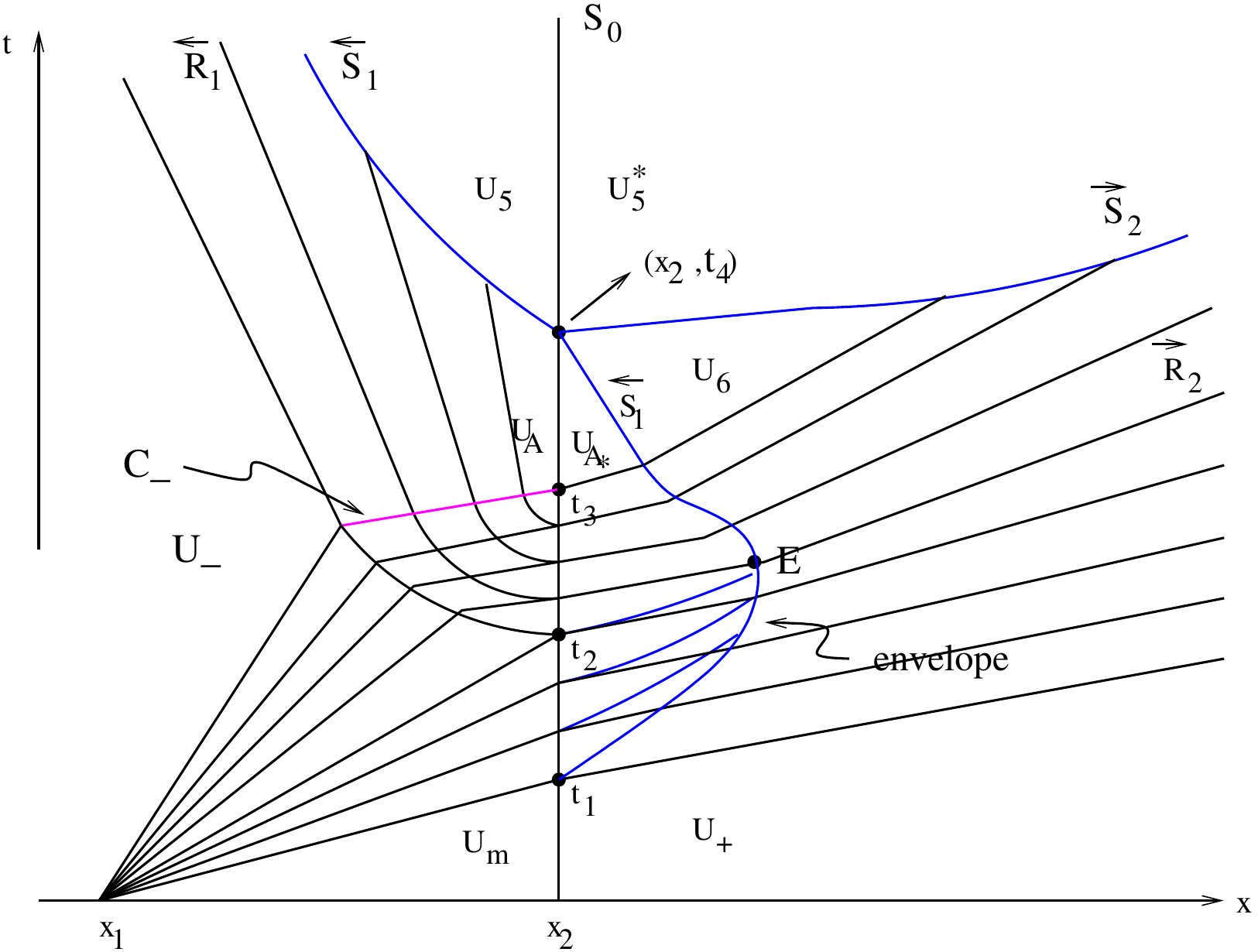}
\end{minipage}
}
\subfigure{
\begin{minipage}[t]{0.47\textwidth}
\centering
\includegraphics[width=0.8725\textwidth]{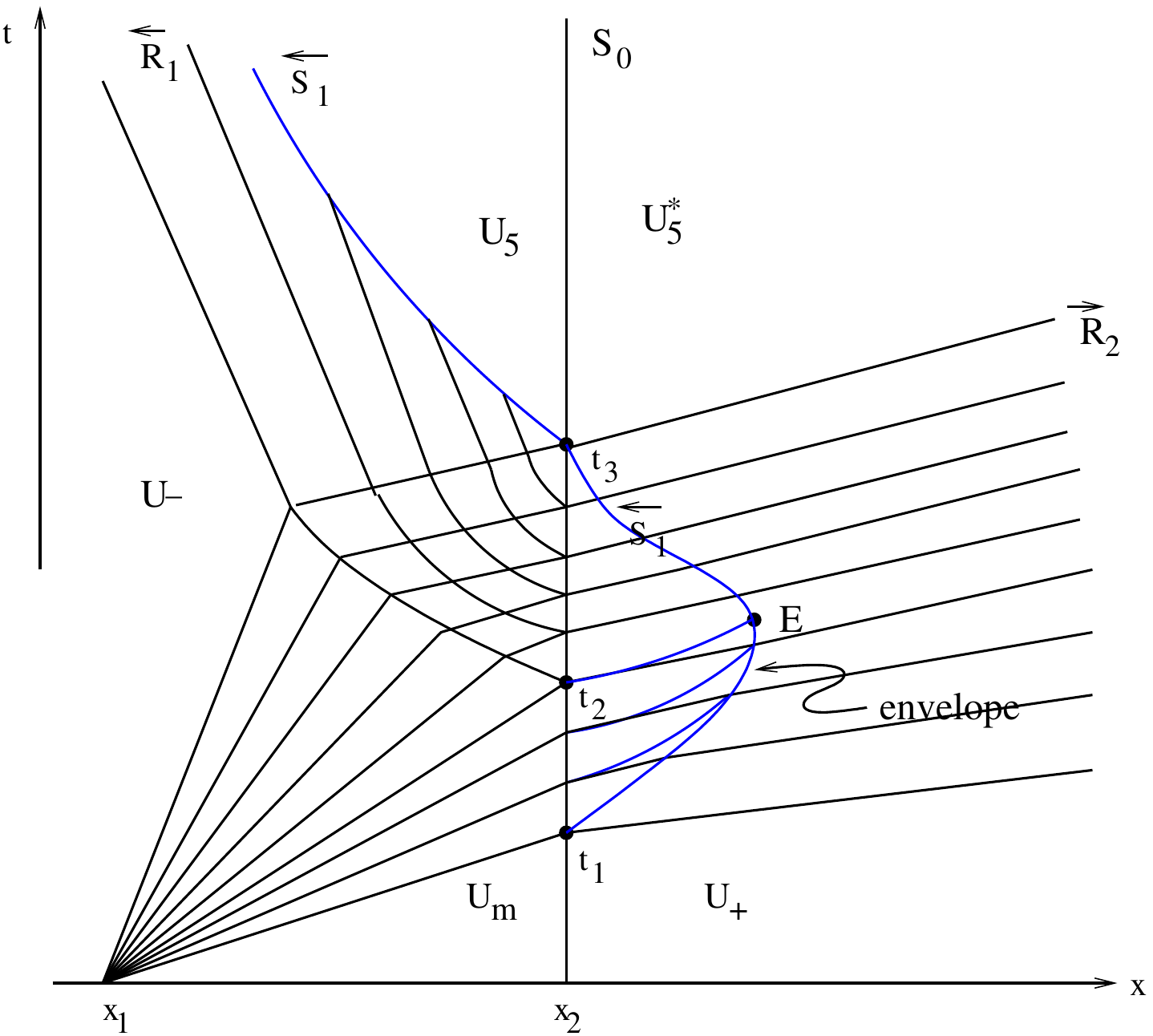}
\end{minipage}
}
\caption*{Fig. 3.11. The results in x-t plane of subcase 3.3(left) and subcase 3.4(right).}
\end{figure}

As $t>t_4$, the refractive shock wave $\overleftarrow{S}_1$ on the left side of $S_0$ begins to penetrate $\overleftarrow{R}_1$ with a varying speed of propagation during the process, i.e., the shock wave curve $\overleftarrow{S}_1: x=x(t)$ is no longer a straight line at $t>t_4$. By using \eqref{2.12} and \eqref{2.18}, the speed of  $\overleftarrow{S}_1$ is determined by the following free boundary problem
\begin{equation}\label{3.25'}
\left\{
\begin{array}{ll}
\displaystyle \frac{{\rm d}x}{{\rm d}t}=u-\sqrt{\frac{g}{2}(h+h_5)\frac{h_5}{h}},\\[6pt]
\displaystyle x-\hat{x}=(u-\sqrt{gh})(t-\hat{t}),\\[6pt]
\displaystyle u+2\sqrt{gh}=u_A+2\sqrt{gh_A},\\[6pt]
\displaystyle x(t_4)=x_2,\quad  h_{A}\leq h \leq h_{-},
\end{array}
\right.
\end{equation}
in which $(\hat{x},\hat{t})$ is on the penetrating characteristic line $C_-$ and representing the translation points from interaction area to $\overleftarrow{R}_1$, see Fig. 3.11. Differentiating the second equation with respect to $t$, one obtains
\begin{equation}\label{3.26'}
\displaystyle \frac{{\rm d}x}{{\rm d}t}=u-\sqrt{gh}+(t-\hat{t})\left(\frac{{\rm d}u}{{\rm d}t}-\frac{1}{2}\sqrt{\frac{g}{h}}\frac{{\rm d}h}{{\rm d}t}\right).
\end{equation}
Combining \eqref{3.26'} with the first equation of \eqref{3.25'}, one has 
\begin{equation}\label{3.27'}
\displaystyle \sqrt{gh}-\sqrt{\frac{g}{2}(h+h_5)\frac{h_5}{h}}=(t-\hat{t})\left(\frac{{\rm d}u}{{\rm d}t}-\frac{1}{2}\sqrt{\frac{g}{h}}\frac{{\rm d}h}{{\rm d}t}\right).
\end{equation}
Besides, from the third equation of \eqref{3.25'}, one has 
\begin{equation}\label{3.28'}
\displaystyle \frac{{\rm d}u}{{\rm d}t}=-\sqrt{\frac{g}{h}}\frac{{\rm d}h}{{\rm d}t}.
\end{equation}
Substituting \eqref{3.28'} into \eqref{3.27'}, one finally obtains
\begin{equation}\label{3.29'}
\displaystyle \frac{{\rm d}h}{{\rm d}t}=-\frac{2}{3(t-\hat{t})}\left(h-\sqrt{\frac{h+h_5}{2}h_5}\right).
\end{equation}
By integrating \eqref{3.29'} on both sides along $x(t)$, it leads to
\begin{equation}\label{3.30'}
\displaystyle {\rm ln}\frac{t-\hat{t}}{t_4-\hat{t}}=\int_{h_A}^{h}\frac{3}{\sqrt{2h_5(h+h_5)}-2h}{\rm d}h.
\end{equation}
It is now clear that $t\rightarrow \infty$ as $h\rightarrow h_5$, therefore $\overleftarrow{S}_1$ cannot penetrate over $\overleftarrow{R}_1$ if $h_5\leq h_-$ holds.
Otherwise, one can expect that $\overleftarrow{S}_1$ penetrates the whole of $\overleftarrow{R}_1$ at the finite time
\begin{equation}\label{3.31'}
\displaystyle t_5=\hat{t}+(t_4-\hat{t}){\rm exp}\left(\int_{h_A}^{h}\frac{3}{\sqrt{2h_5(h+h_5)}-2h}{\rm d}h\right).
\end{equation}
For completeness, we conclude that the refractive shock wave $\overleftarrow{S}_1$ on the left side of $S_0$ is able to cross the whole of 
$\overleftarrow{R}_1$ for $h_5>h_-$, whereas it cannot for $h_5\leq h_-$ and eventually has $x-x_2=(u_5-\sqrt{gh_5})(t-t_4)$ as its asymptote.

Now we are left with the interaction of $\overrightarrow{S}_2$ and $\overrightarrow{R}_2$ on the right side of $S_0$ as $t>t_4$. Actually,  the speed of  $\overrightarrow{S}_2$ in this case is determined by
\begin{equation}\label{3.32'}
\left\{
\begin{array}{ll}
\displaystyle \frac{{\rm d}x}{{\rm d}t}=u+\sqrt{\frac{g}{2}(h+h_{5*})\frac{h_{5*}}{h}},\\[6pt]
\displaystyle x-\tilde{x}=(u+\sqrt{gh})(t-\tilde{t}),\\[6pt]
\displaystyle u-2\sqrt{gh}=u_6-2\sqrt{gh_6},\\[6pt]
\displaystyle x(t_4)=x_2,\quad  h_{6}\leq h \leq h_{+},
\end{array}
\right.
\end{equation}
By applying the same process of \eqref{3.26'}-\eqref{3.30'} to \eqref{3.32'}, we conclude that the shock wave $\overrightarrow{S}_2$ is able to cross the whole of 
$\overrightarrow{R}_2$ for $h_{5*}>h_+$, whereas it cannot for $h_{5*}\leq h_{+}$ and eventually has $x-x_2=(u_{5*}+\sqrt{gh_{5*}})(t-t_4)$ as its asymptote.

In order to straighten out the whole structure, we state that when the shock waves penetrate the whole of rarefaction waves, the solution in the large time scale includes a  backward shock wave $\overleftarrow{S}_1$ from $U_-$ and $U_5$, followed by a stationary wave from $U_5$ to $U_{5}^{*}$, finally followed by a forward shock wave $\overrightarrow{S}_2(U_+,U_5^{*})$.

\noindent
{\bf Subcase 3.4.} This case happens when $U_6$ and $U_{5*}$ coincide, the 1-shock wave $\overleftarrow{S}_1$ in subcase 3.3 will leave inside $\overrightarrow{R}_2$ before it touches $S_0$. There will be no forward shock wave reflecting in this case. On the left side of $S_0$, the refractive shock wave $\overleftarrow{S}_1$ will interact with $\overleftarrow{R}_1$ as $t>t_3$. By a similar discussion as \eqref{3.25'} to \eqref{3.31'}, one can conclude here that $\overleftarrow{S}_1$ will either overtake the whole of $\overleftarrow{R}_1$ at a finite time or eventually leave inside $\overleftarrow{R}_1$ and  have a curve as its asymptote. The result is shown in Fig. 3.11 (right).

\begin{figure}[h]
\subfigure{
\begin{minipage}[t]{0.45\textwidth}
\centering
\includegraphics[width=0.95\textwidth]{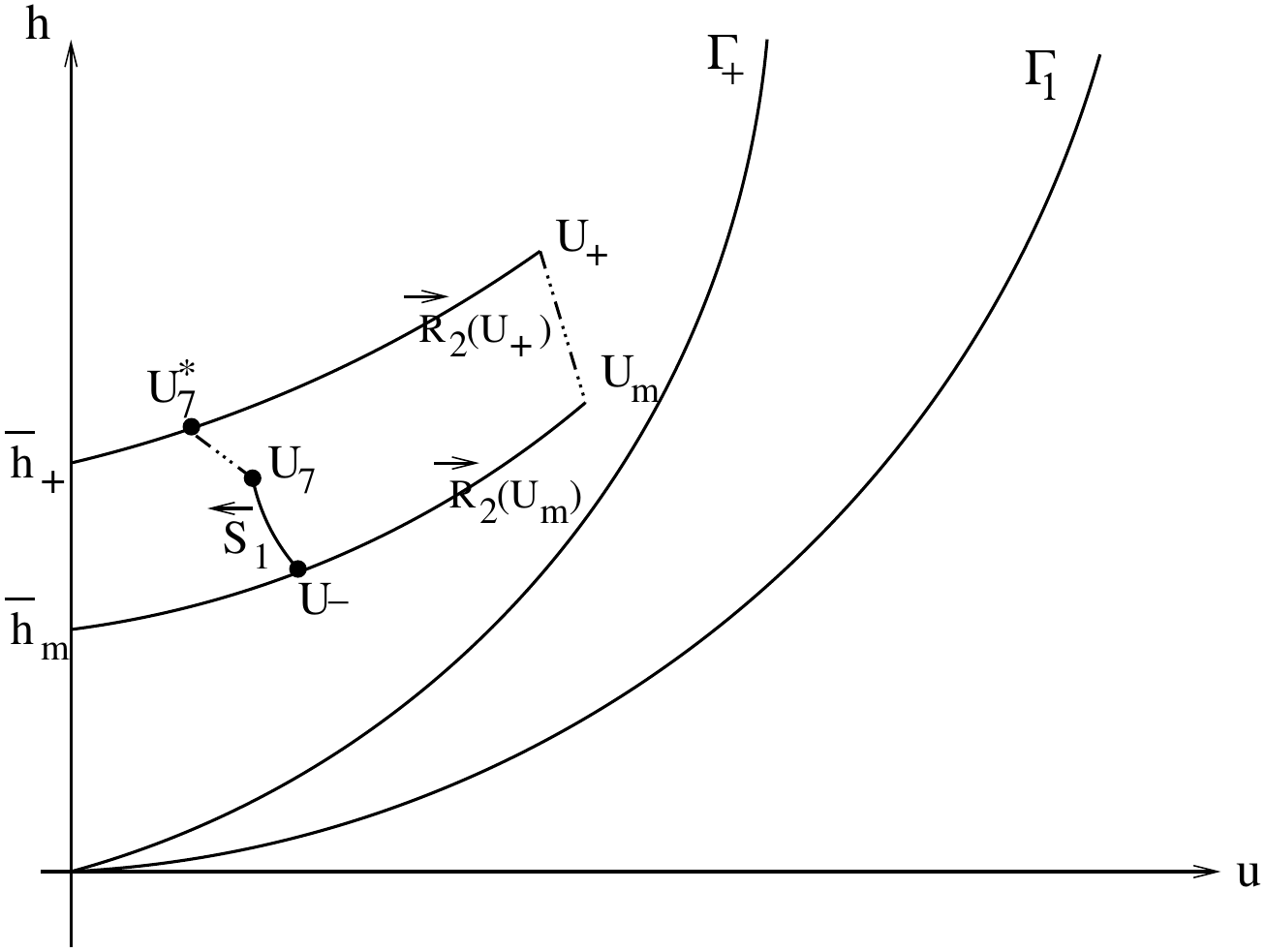}
\end{minipage}
}
\subfigure{
\begin{minipage}[t]{0.45\textwidth}
\centering
\includegraphics[width=0.95\textwidth]{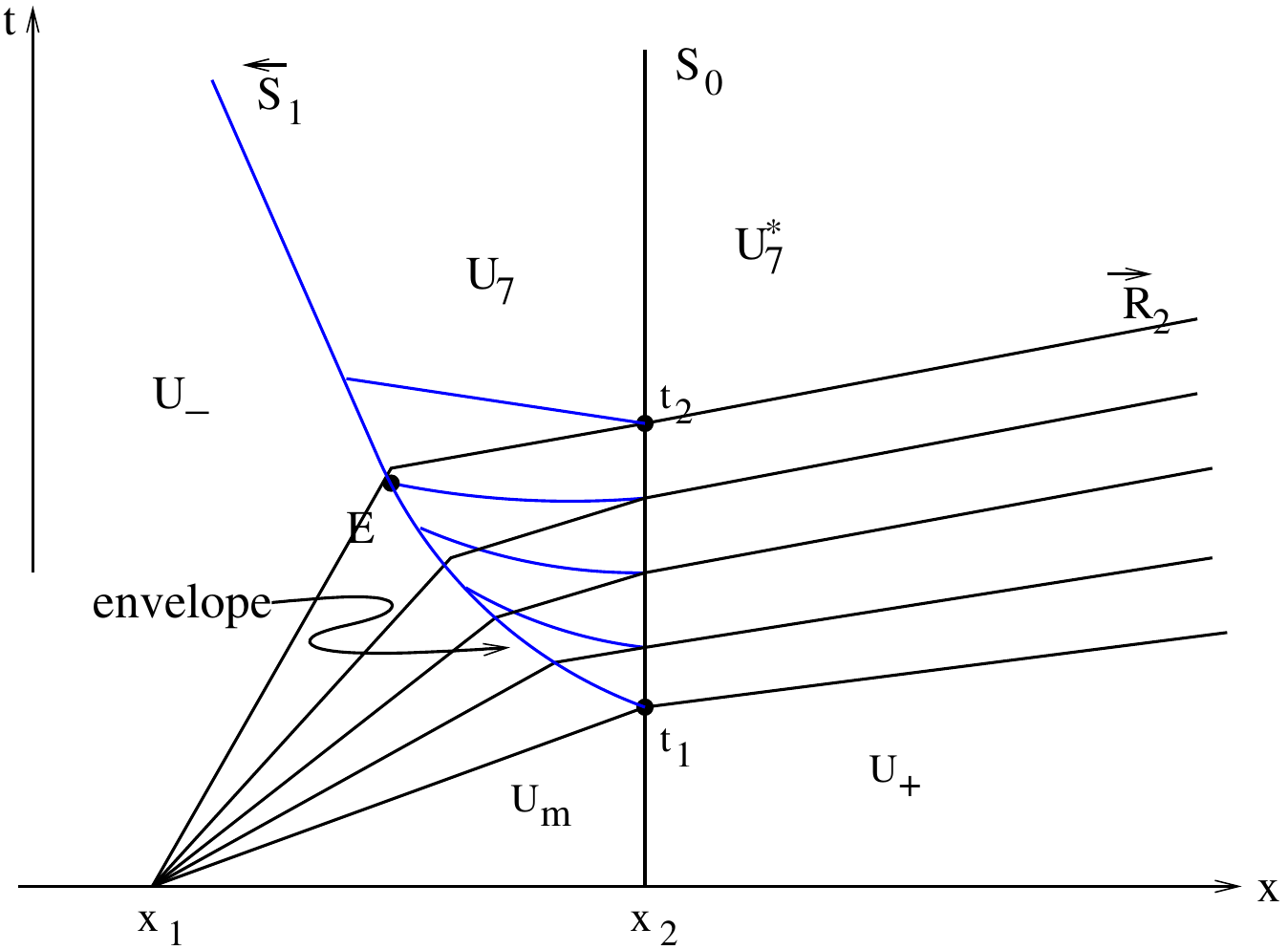}
\end{minipage}
}
\caption*{Fig. 3.12. Case 4. $U_m\in {\rm I\!I\!I}$, $U_+\in {\rm I\!I\!I}$.}
\end{figure}

\noindent
{\bf Case 4.} $U_m\in {\rm I\!I\!I}$, $U_+\in {\rm I\!I\!I}$. See Fig. 3.12.\\
In the subcritical area, $S_0(U_1,U_0)$ is below $\overrightarrow{R}_2(U_+,U)$, touches $h-$axis at $(0,\tilde{h})$, see Fig. 3.6.

The interaction process is from $t_1$ to $t_2$. Denote $U_0^{s}\in \overleftarrow{S}_1(U,U_-)$, and $U_7^{*}=S_0(U_0^{s*},U_0^{s})\cap \overrightarrow{R}_2(U_+,U)$, which is obtained by $U_7\in \overleftarrow{S}_1(U,U_-)$. We solve the initial boundary value problem \eqref{2.11} with
\begin{equation}\label{3.33}
(u,h)=\left\{
\begin{array}{l}
\displaystyle \big(\xi+\sqrt{gh}, \big(\frac{u_-+2\sqrt{gh_-}-\xi}{3\sqrt{g}}\big)^2\big),\quad (x,t)\in C_-,~u_--\sqrt{gh_-}<\xi<u_m-\sqrt{gh_m}, \\[8pt]
\displaystyle U_0=(u_0,h_0), x=x_2,\quad t_2<t<t_3,\\[7pt]
\displaystyle (u_1,h_1),\quad x=x_2,t>t_3,
\end{array}
\right.
\end{equation}
where $U_0\in S_0(U_1,U_0)$ is obtained by $U_1\in \overrightarrow{R}_2(U_+,U)$ ($U_1$ is the right-hand state). The solution of \eqref{3.33} contains a reflecting compressible rarefaction wave, which will develop a shock wave propagates on the left. See Fig. 3.12 (right). The large time behavior of the solution from the $p-$system theory is
\begin{equation}
U_-\oplus \overleftarrow{S}_1(U_1,U_-)\oplus S_0(U_{1}^{*},U_1)\oplus  \overrightarrow{R}_2(U_+,U_1^{*})\oplus U_+.
\end{equation}

\subsection{Interaction of a shock wave with a stationary wave}
In this section, we consider the interaction of shock wave with stationary wave. In the initial value problem \eqref{3.1}, we have $U_m\in \overrightarrow{S}_2(U,U_-)$, $U_+\in S_0(U,U_m)$, i.e,
\begin{equation}\label{3.34}
\left\{\begin{array}{ll}
\displaystyle U_m\in \overrightarrow{S}_2(U,U_-):~ u=u_--\sqrt{\frac{g}2(h-h_-)^2\left(\frac1{h}+\frac1{h_-}\right)},\quad h<h_-.\\[5pt]
W_3(U_+;U_m):~\left\{\begin{array}{lll}
h_m u_m = h_+u_+,\\[5pt]
\displaystyle \frac{u_m^2}2+g(h_m+a_0)=\frac{u_+^2}2+g(h_++a_1), &a_1<a_{max}.
\end{array}
\right.
\end{array}\right.
\end{equation}
The state $U=(u,h)\in \overrightarrow{S}_2(U_+,U)$($U_+$ is the right-hand state) is given by
\begin{equation}
u=u_++\sqrt{\frac{g}2(h-h_+)^2\left(\frac1{h}+\frac1{h_+}\right)},\quad h>h_+.
\end{equation}
Liu(\cite{Liu2}) has proved that shock wave tends to decelerate for flows along an expanding duct. Since shallow water equations \eqref{1.1} have similar structures with the isentropic Euler equations, we apply the conclusion to the shallow water equations directly, i.e., the shock wave decelerates along the direction of decreasing $h$. Next, we discuss the interaction results case by case.

\noindent
{\bf Case 1.} $U_m, U_+$ are on the right side of $\Gamma_+$. We have $u_m>\sqrt{gh_m}$, $u_+>\sqrt{gh_+}$ in this case. When $\overrightarrow{S}_2(U_m,U_-)$ overtakes $S_0(U_+,U_m)$, one solves a new Riemann problem of \eqref{1.1} with
\begin{equation}\label{3.35}
(u,h,a)\big|_{t=1}=
\left\{\begin{array}{ll}
U_-=(u_-,h_-,a_0),   x<x_2,\\[5pt]
U_+=(u_+,h_+,a_1),   x>x_2.
\end{array}\right.
\end{equation}
In the supercritical area, the left-hand state $U_-$ will first jump to $U_{-*}$ by $S_0$. According to the relative positions of $U_{-*}$ and $\overrightarrow{S}_2(U_+,U)$, we discuss as follows.

\begin{figure}[h]
\subfigure{
\begin{minipage}[t]{0.45\textwidth}
\centering
\includegraphics[width=0.9\textwidth]{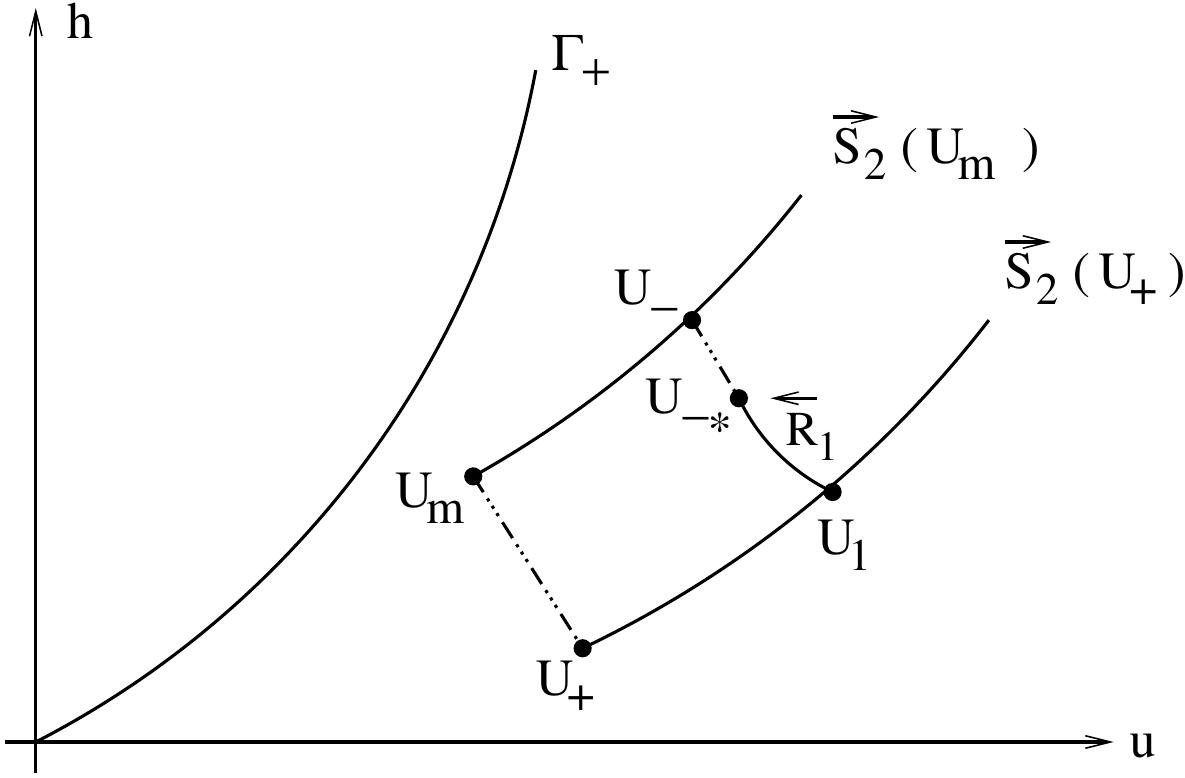}
\end{minipage}
}
\subfigure{
\begin{minipage}[t]{0.45\textwidth}
\centering
\includegraphics[width=0.9\textwidth]{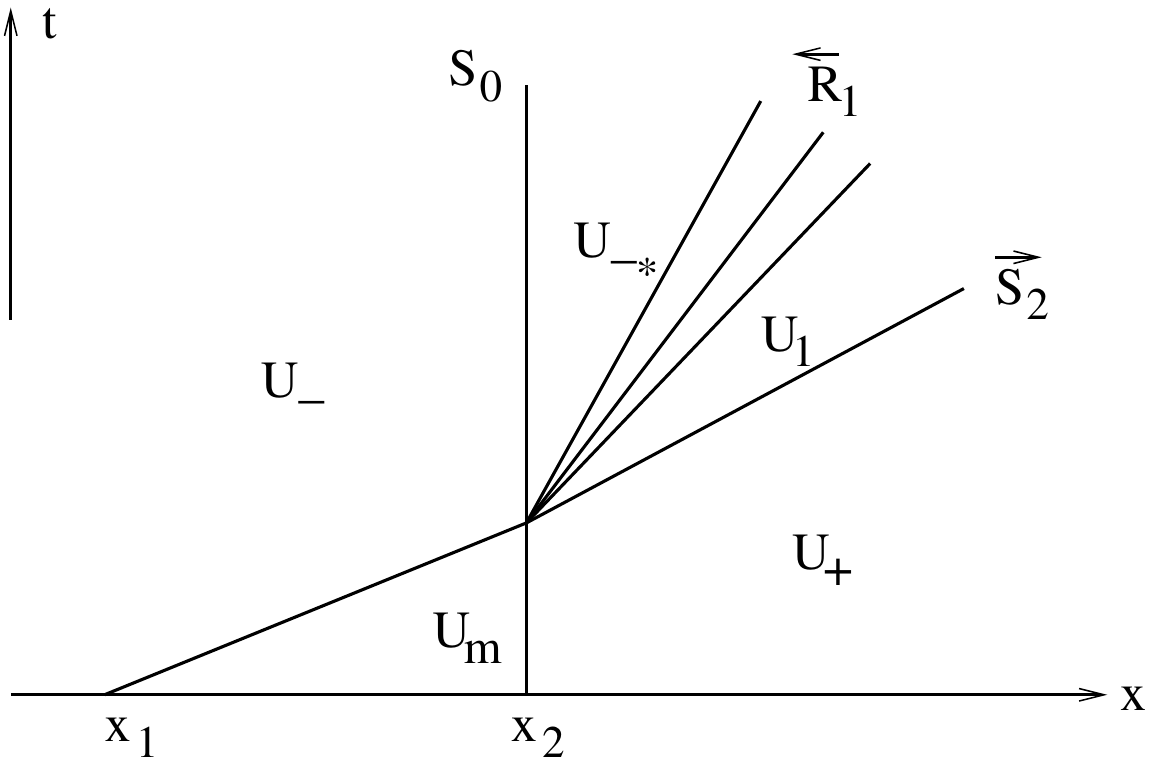}
\end{minipage}
}
\caption*{Fig. 3.13. Subcase 1.1. $U_{-*}$ is on the left of $\overrightarrow{S}_2(U_+,U)$.}
\end{figure}

\noindent
{\bf Subcase 1.1.} $U_{-*}$ is on the left of $\overrightarrow{S}_2(U_+,U)$. In this case $U_-$ jumps to $U_{-*}$ by $S_0$ first, then $U_{-*}$ connects with $U_1$ by $\overleftarrow{R}_1(U,U_{-*})$, finally $U_1$ jumps to $U_+$ by $\overrightarrow{S}_2(U_+,U_1)$. See Fig. 3.13. That is 
\begin{equation}
\overrightarrow{S}_2(U_m,U_-)\oplus S_0(U_+,U_m)\rightarrow S_0(U_{-*},U_-) \oplus \overleftarrow{R}_1(U_1,U_{-*})  \oplus \overrightarrow{S}_2(U_+,U_1),
\end{equation}
which means the forward shock wave will transmit a backward rarefaction wave when it penetrates the stationary wave.

\noindent
{\bf Subcase 1.2.} $U_{-*}$ is on the right of $\overrightarrow{S}_2(U_+,U)$. Then $U_-$ jumps to $U_{-*}$ by $S_0$ first, $U_{-*}$ connects with $U_2$ by $\overleftarrow{S}_1(U,U_{-*})$, finally $U_2$ jumps to $U_+$ by $\overrightarrow{S}_2(U_+,U_2)$. See Fig. 3.14. That is
\begin{equation}
\overrightarrow{S}_2(U_m,U_-)\oplus S_0(U_+,U_m)\rightarrow S_0(U_{-*},U_-)\oplus \overleftarrow{S}_1(U_2,U_{-*})\oplus \overrightarrow{S}_2(U_+,U_2),
\end{equation}
which means the forward shock wave will transmit a backward shock wave when it penetrates the stationary wave.

\begin{figure}
\subfigure{
\begin{minipage}[t]{0.45\textwidth}
\centering
\includegraphics[width=0.9\textwidth]{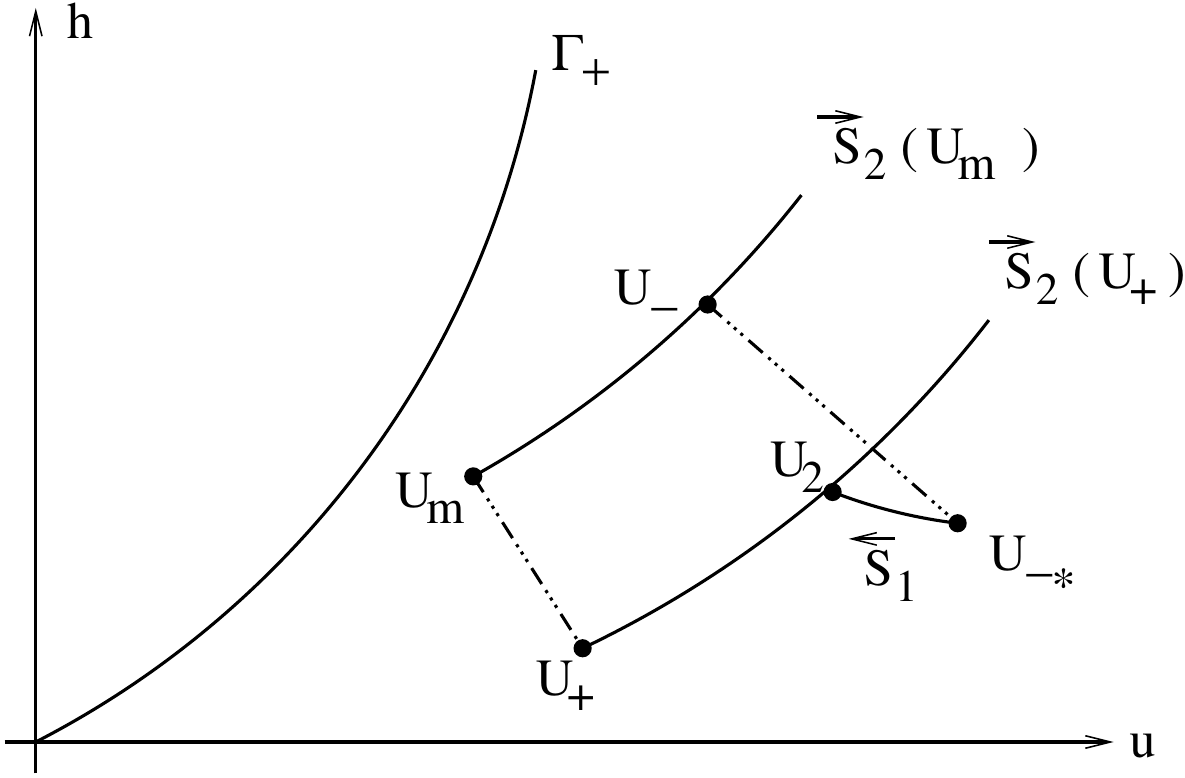}
\end{minipage}
}
\subfigure{
\begin{minipage}[t]{0.45\textwidth}
\centering
\includegraphics[width=0.9\textwidth]{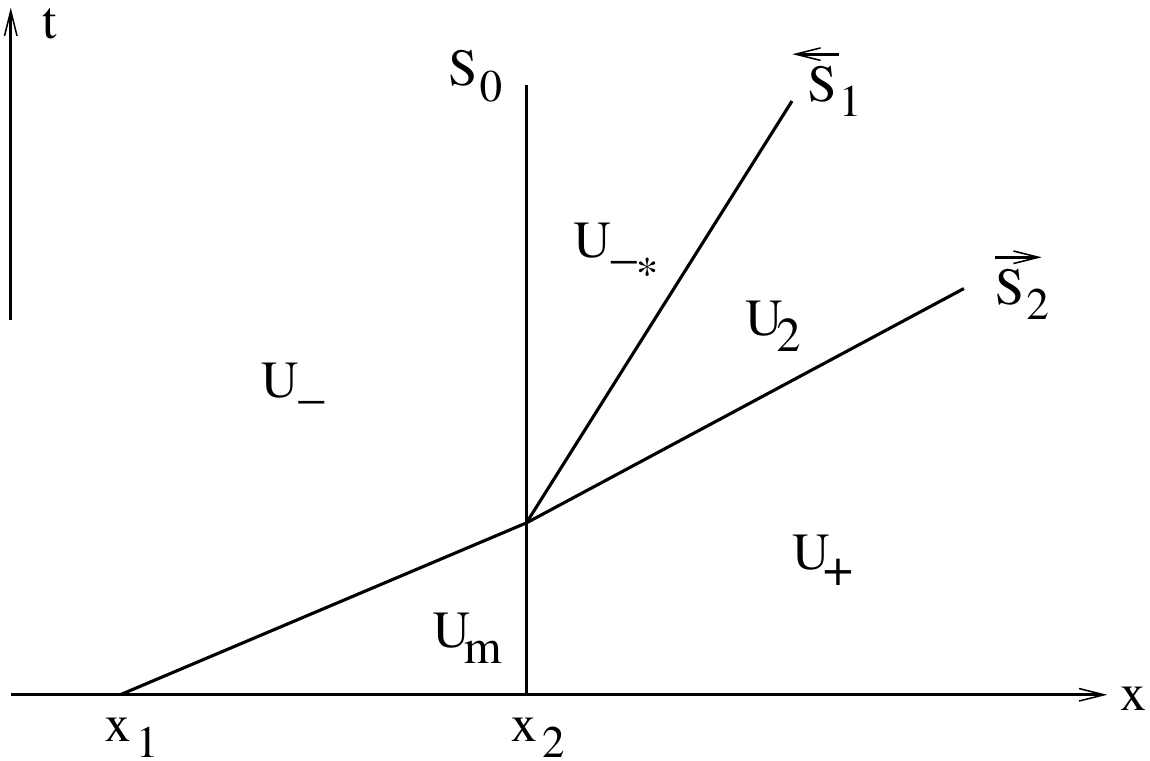}
\end{minipage}
}
\caption*{Fig. 3.14. Subcase 1.2. $U_{-*}$ is on the right of $\overrightarrow{S}_2(U_+,U)$.}
\end{figure}

\noindent
{\bf Case 2.} $U_m\in D_2, U_+\in D_2$ and $U_-\in D_2$, where $D_2$ is defined in \eqref{2.7'}. We choose the subcritical solution when using the stationary wave. The results are discussed as follows.

\noindent
{\bf Subcase 2.1.} $U_{-*}$ is below $\overrightarrow{S}_2(U_+,U)$. Then $U_-$ jumps to $U_3$ by $\overleftarrow{S}_1(U_3,U_-)$, $U_3$ connects with $U_3^{*}$ by $S_0(U_3^{*},U_3)$, and finally $U_3^{*}$ jumps to $U_+$ by $\overrightarrow{S}_3(U_+,U_3^{*})$. See Fig. 3.15. That is
\begin{equation}
\overrightarrow{S}_2(U_m,U_-)\oplus S_0(U_+,U_m)\rightarrow \overleftarrow{S}_1(U_3,U_-) \oplus S_0(U_3^{*},U_3)  \oplus \overrightarrow{S}_2(U_+,U_3^{*}),
\end{equation}
which means the forward shock wave will reflect a backward shock wave when it penetrates the stationary wave.

\noindent
{\bf Subcase 2.2.} $U_{-*}$ is above $\overrightarrow{S}_2(U_+,U)$. Then $U_-$ passes a backward rarefaction wave first. Denote $U_c=\Gamma_+\cap \overleftarrow{R}_1(U,U_-)$. See Fig. 3.16.

If $U_{c}^{*}$ is below $\overrightarrow{S}_2(U_+,U)$, then $U_-$ connects with $U_4\in D_1$ by $\overleftarrow{R}_1(U,U_-)$, $U_4$ jumps to $U_4^{*}\in \overrightarrow{S}_2(U_+,U)$ by a stationary wave $S_0(U_4^{*},U_4)$, finally $U_4^{*}$ jumps to $U_+$ by $\overrightarrow{S}_2(U_+,U_4^{*})$. See Fig. 3.16. That is
\begin{equation}
\overrightarrow{S}_2(U_m,U_-)\oplus S_0(U_+,U_m)\rightarrow \overleftarrow{S}_1(U_4,U_-) \oplus S_0(U_4^{*},U_4)  \oplus \overrightarrow{S}_2(U_+,U_4^{*}),
\end{equation}
which means the forward shock wave will reflect a backward rarefaction wave when it penetrates the stationary wave.

\begin{figure}
\subfigure{
\begin{minipage}[t]{0.45\textwidth}
\centering
\includegraphics[width=0.9\textwidth]{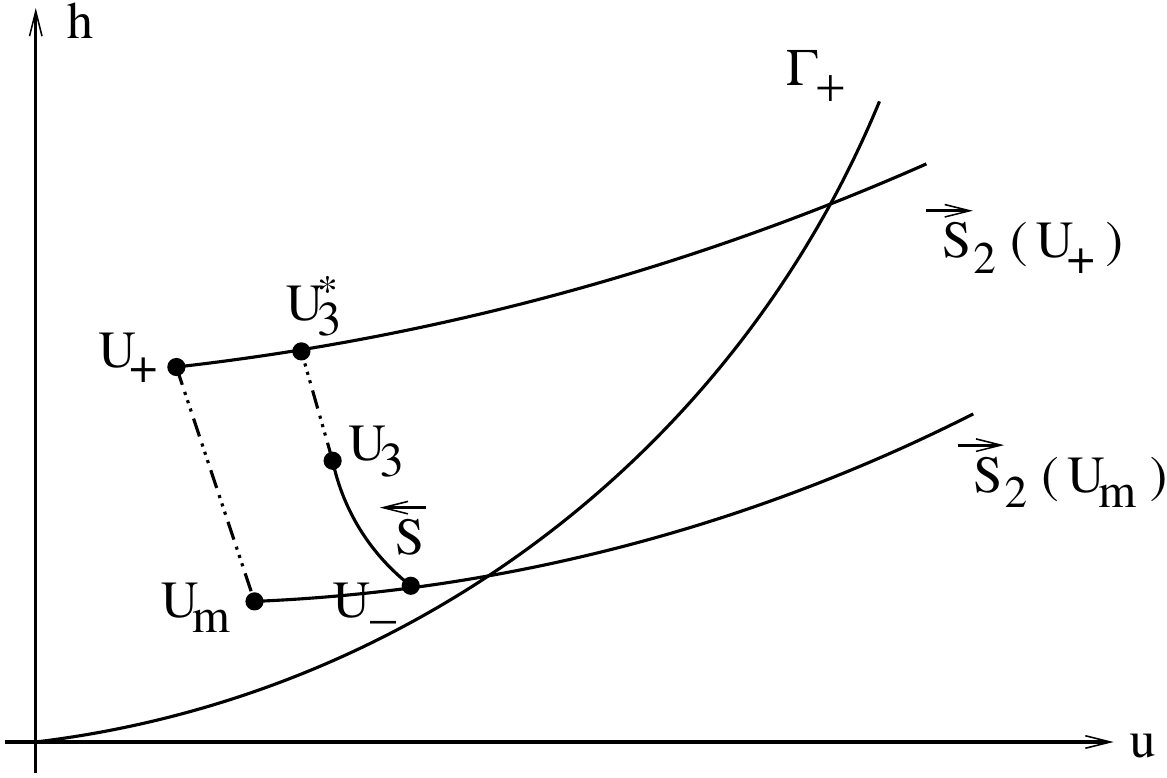}
\end{minipage}
}
\subfigure{
\begin{minipage}[t]{0.45\textwidth}
\centering
\includegraphics[width=0.9\textwidth]{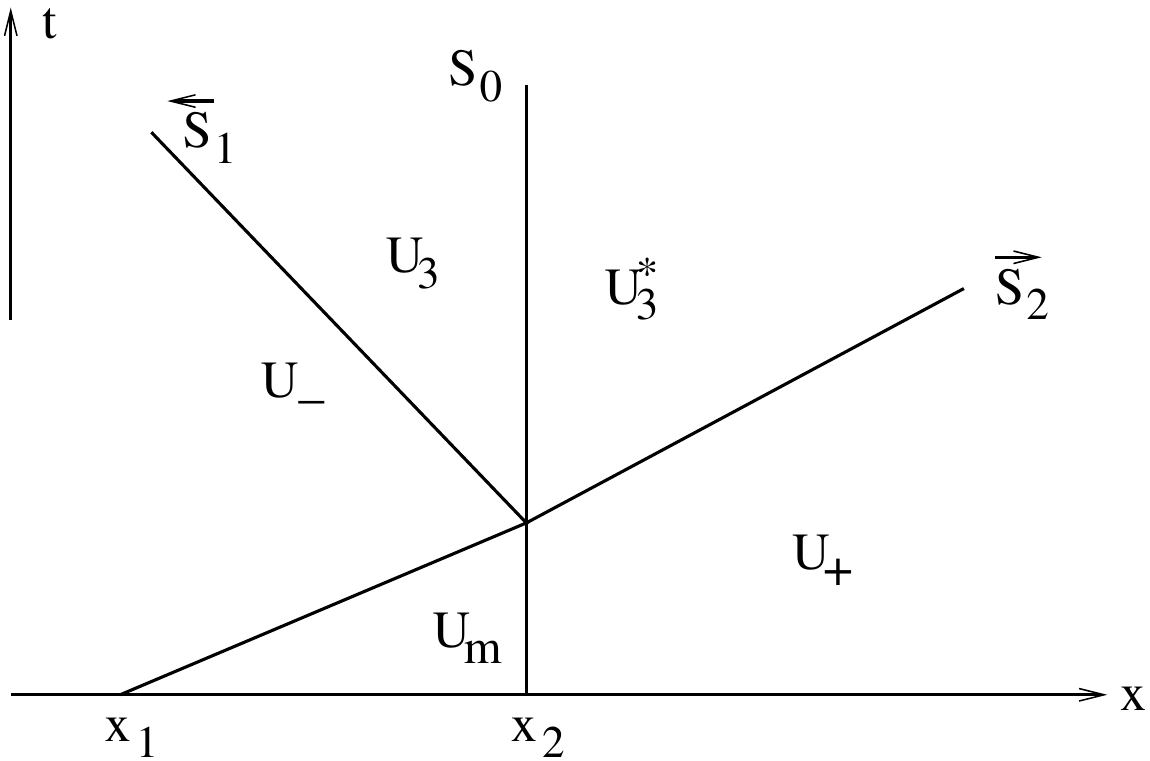}
\end{minipage}
}
\caption*{Fig. 3.15. Subcase 2.1. $U_{-*}$ is below $\overrightarrow{S}_2(U_+,U)$.}
\end{figure}

\begin{figure}
\subfigure{
\begin{minipage}[t]{0.45\textwidth}
\centering
\includegraphics[width=0.9\textwidth]{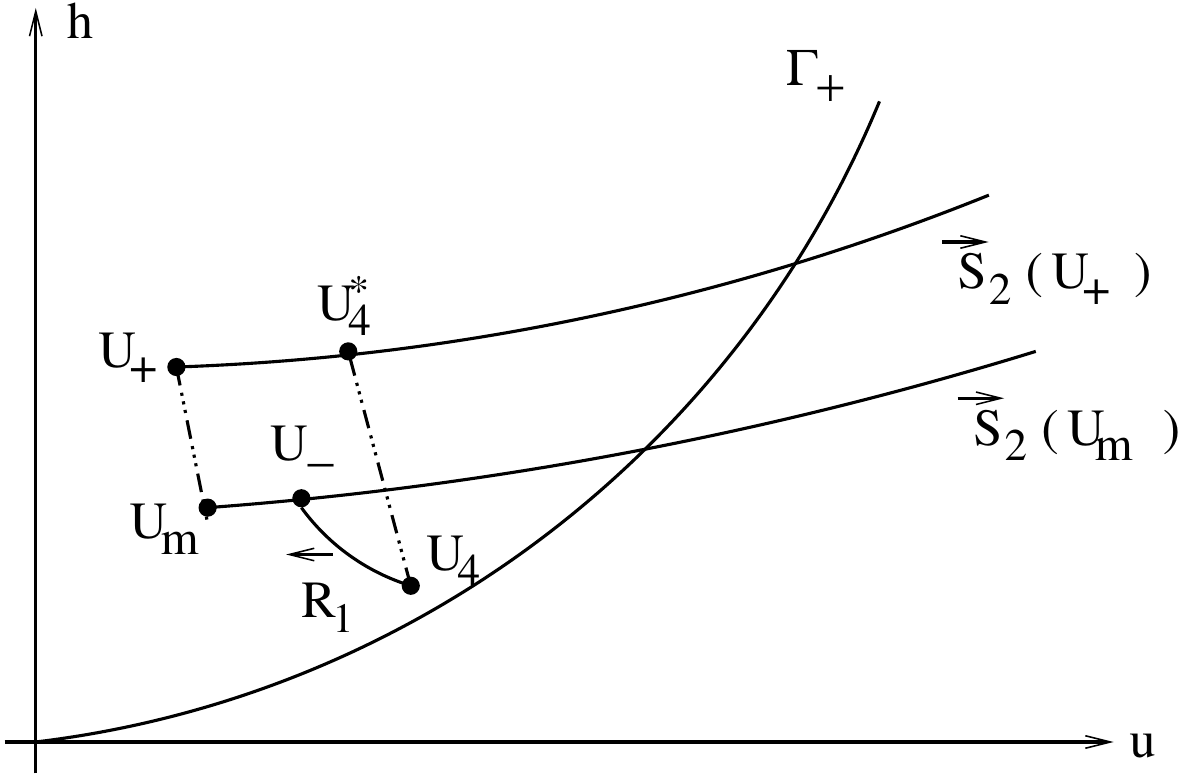}
\end{minipage}
}
\subfigure{
\begin{minipage}[t]{0.45\textwidth}
\centering
\includegraphics[width=0.9\textwidth]{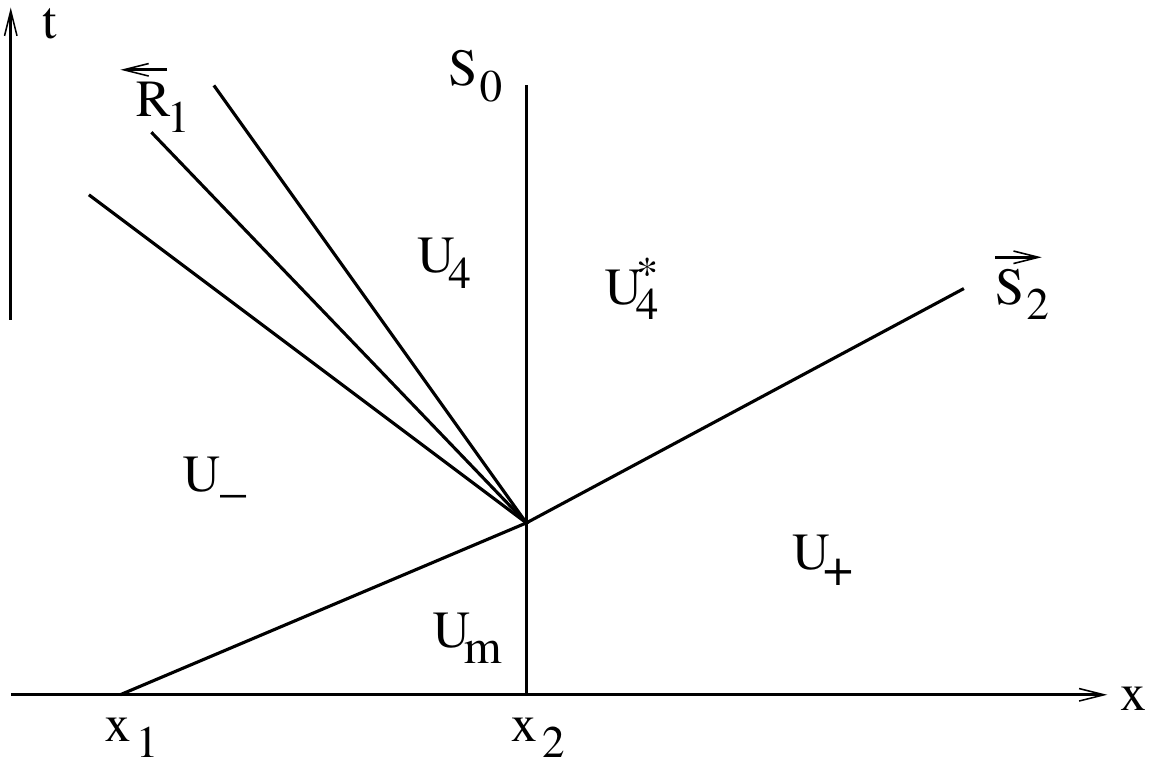}
\end{minipage}
}
\caption*{Fig. 3.16. Subcase 2.2. $U_{-*}$ is above $\overrightarrow{S}_2(U_+,U)$.}
\end{figure}

\noindent
{\bf Subcase 2.3.}
If $U_c^{*}$ is above $\overrightarrow{S}_2(U_+,U)$. Then $U_-$ connects with $U_c$ by
$\overleftarrow{R}_1(U,U_-)$ first, $U_c$ jumps to $U_{c*}$ by $S_0(U_{c*},U_c)$, then $U_{c*}$ connects with $U_5$ by $\overleftarrow{R}_1(U_5,U_{c*})$, finally $U_5$ jumps to $U_+$ by $\overrightarrow{S}_2(U_+,U_5)$. See Fig. 3.17.
This case means the forward shock wave will reflect a backward rarefaction wave which coincides with the stationary wave $S_0$.

\begin{figure}[h]
\subfigure{
\begin{minipage}[t]{0.45\textwidth}
\centering
\includegraphics[width=0.9\textwidth]{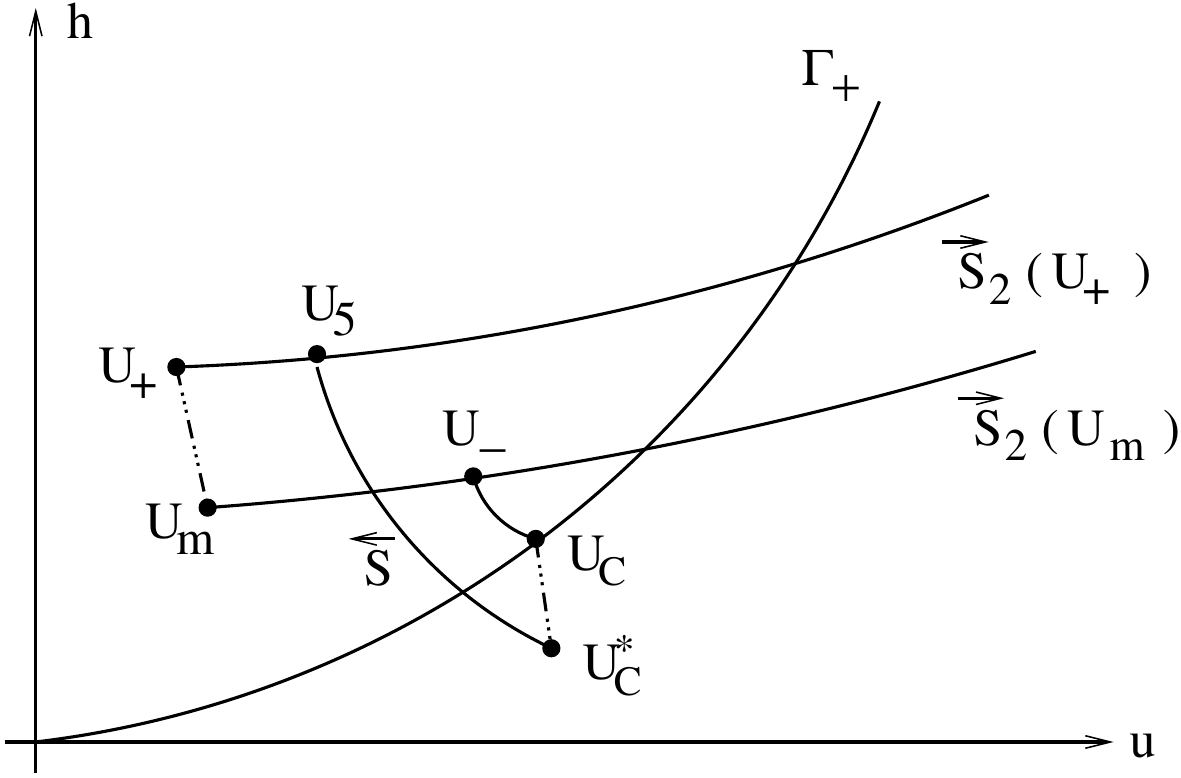}
\end{minipage}
}
\subfigure{
\begin{minipage}[t]{0.45\textwidth}
\centering
\includegraphics[width=0.9\textwidth]{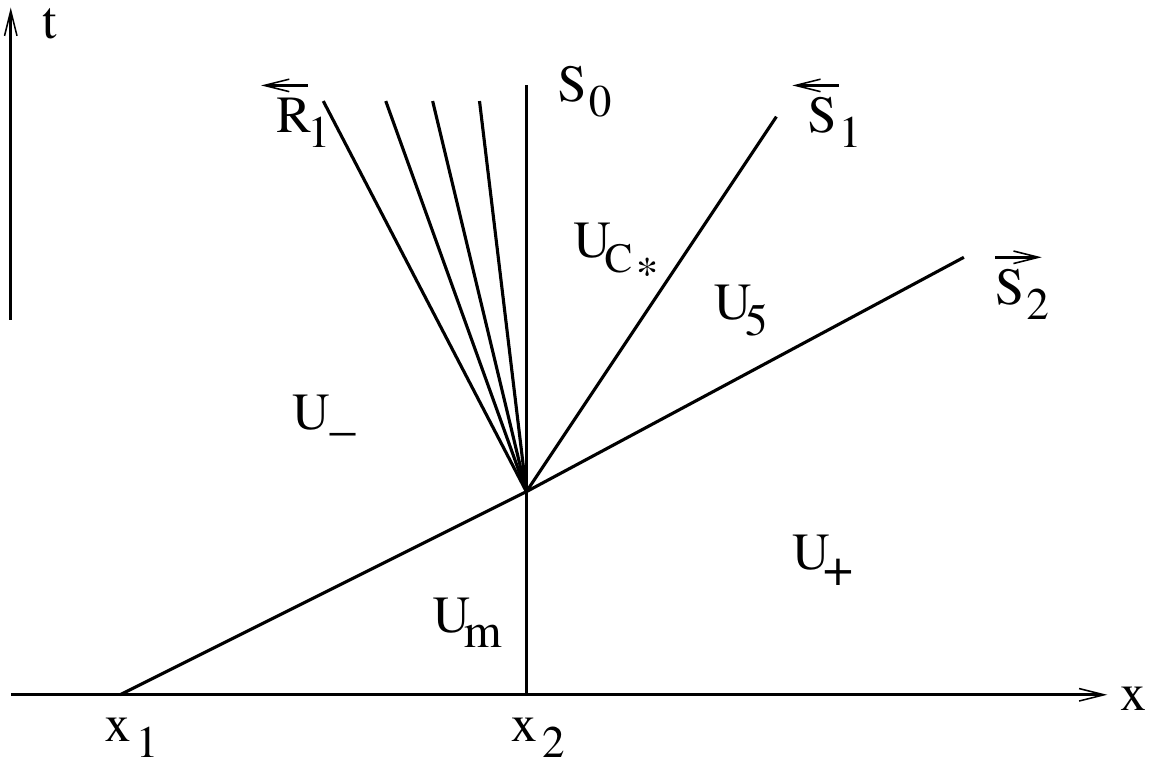}
\end{minipage}
}
\caption*{Fig. 3.17. Subcase 2.3. $U_{c}^{*}$ is above $\overrightarrow{S}_2(U_+,U)$.}
\end{figure}

\noindent
{\bf Case 3.} $U_m\in D_2$, $U_+\in D_2$ and $U_-\in D_3$. We discuss this case in the following two subcases.

\noindent
{\bf Subcase 3.1.} Denote $\widetilde{U}_-\in \overleftarrow{S}_1(U,U_-)$ which satisfies $\sigma_1(\widetilde{U}_-,U_-)=0$, $\widetilde{U}_-$ jumps to $\widetilde{U}_-^{*}$ by a stationary wave $S_0$. If $\widetilde{U}_-^{*}$ is below $\overrightarrow{S}_3(U_+,U)$, then solution is:
 $U_-$ connects with $U_6$ by $\overleftarrow{S}_1$ first, $U_6$ jumps to
 $U_6^{*}\in\overrightarrow{S}_2(U_+,U)$ by a stationary wave, finally $U_6^{*}$ connects with $U_+$ by $\overrightarrow{S}_2(U_+,U)$. See Fig. 3.18. 
\begin{equation}
\overrightarrow{S}_2(U_m,U_-)\oplus S_0(U_+,U_m)\rightarrow \overleftarrow{S}_1(U_6,U_-) \oplus S_0(U_6^{*},U_6)  \oplus \overrightarrow{S}_2(U_+,U_6^{*}).
\end{equation}

\begin{figure}
\subfigure{
\begin{minipage}[t]{0.45\textwidth}
\centering
\includegraphics[width=0.9\textwidth]{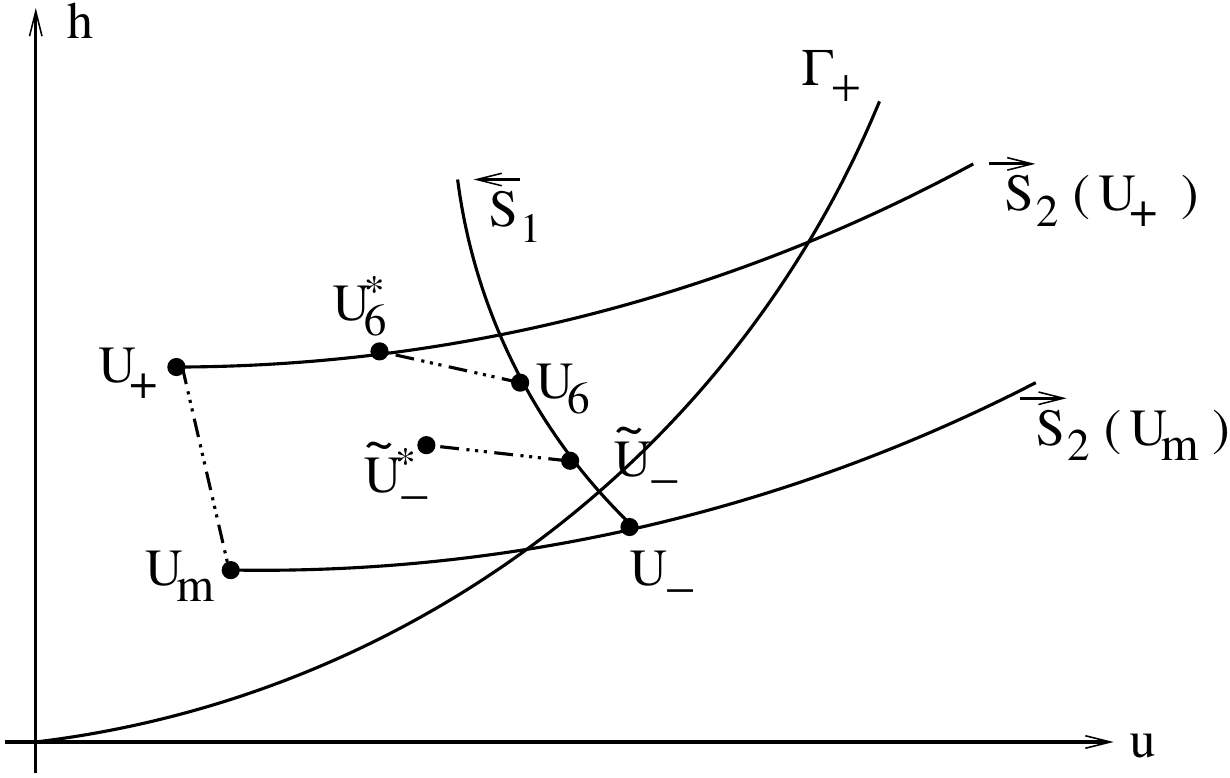}
\end{minipage}
}
\subfigure{
\begin{minipage}[t]{0.45\textwidth}
\centering
\includegraphics[width=0.9\textwidth]{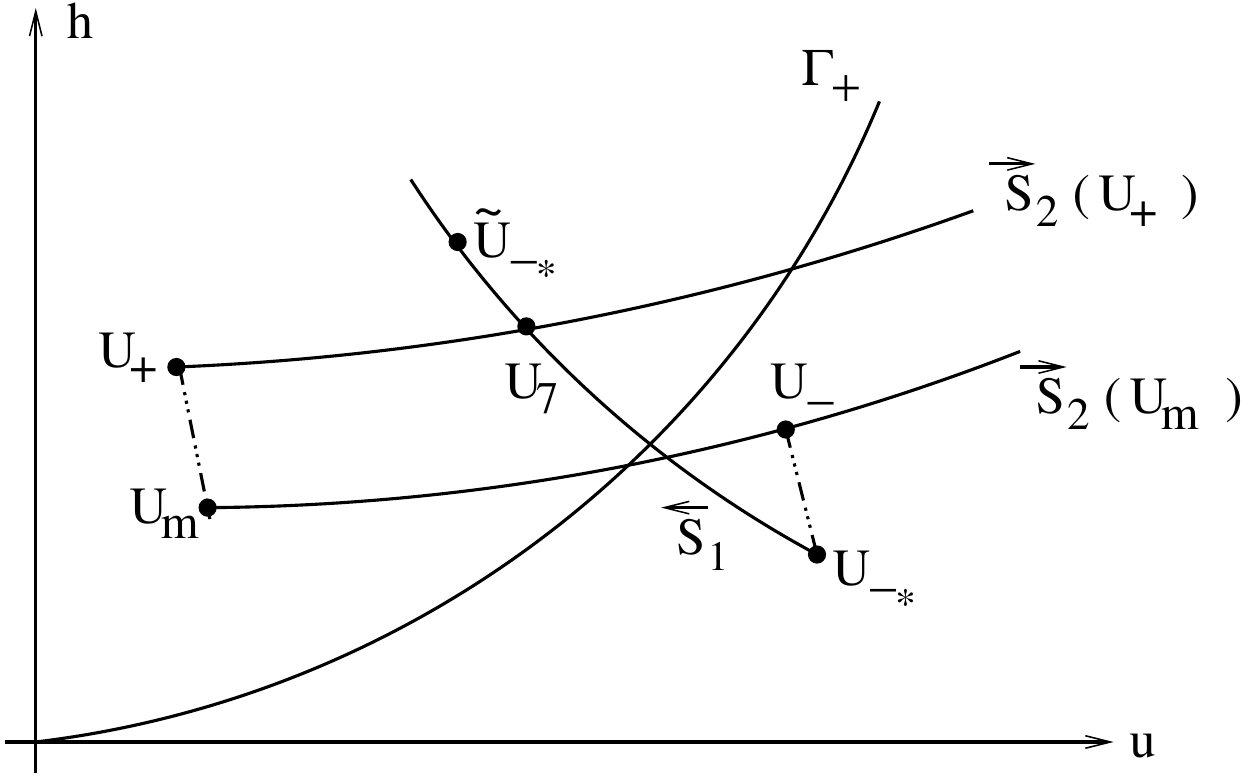}
\end{minipage}
}
\caption*{Fig. 3.18. Subcase 3.1. $U_{-*}$ is below $\overrightarrow{S}_2(U_+,U)$.}
\end{figure}

\noindent
{\bf Subcase 3.2.} $U_-$ jumps to $U_{-*}$ by $S_0(U_{-*},U_-)$. There exists a point $\widetilde{U}_{-*}\in D_1$ which satisfies $\sigma_1(\widetilde{U}_{-*},U_{-*})=0$. If $\widetilde{U}_{-*}$ is above $\overrightarrow{S}_2(U_+,U)$, then the solution is: $U_-$ jumps to $U_{-*}$ by a stationary wave $S_0$ first, $U_{-*}$ connects with $U_7\in \overrightarrow{S}_2(U_+,U)$ by $\overleftarrow{S}_1(U_7,U_{-*})$, finally $U_7$ jumps to $U_+$ by $\overrightarrow{S}_2(U_+,U)$. This case is similar with subcase 1.2. That is
\begin{equation}
\overrightarrow{S}_2(U_m,U_-)\oplus S_0(U_+,U_m)\rightarrow S_0(U_{-*},U_-) \oplus \overleftarrow{S}_1(U_7,U_{-*})  \oplus \overrightarrow{S}_2(U_+,U_7).
\end{equation}

In summary, we have obtained the results of interaction of stationary wave with shock waves or rarefaction waves for the shallow water equations. When a forward rarefaction wave interacts with a stationary wave, it will transmit a forward rarefaction wave. At the same time, a backward wave will transmit or reflect as well. We further discuss a more complicated case in which a compressible backward wave is transmitted. By solving free boundary problems, the results in large time scales are well investigated. When a forward shock wave interacts with a stationary wave, it will penetrate the stationary wave, and either transmit or reflect a backward wave. 

\vskip2mm
{\bf Acknowledgments}.
This work is partially supported by NSFC 11771274. The authors are very
grateful to the anonymous referees for their corrections and suggestions, which improved the original manuscript greatly.


\vskip2mm

 
\begin{thebibliography}{20}

\bibitem{Alcrudo} F. Alcrudo, F. Benkhaldoun,  {\em Exact solutions to the Riemann problem of the shallow water equations with a bottom step}, Comput. Fluids, 30: 643-671, 2001.

\bibitem{Andrianov1} N. Andrianov, G. Warnecke, {\em On the solution to the Riemann problem for the compressible duct flow}, SIAM J. Appl. Math., 64: 878-901, 2004.

\bibitem{Andrianov2} N. Andrianov, G. Warnecke, {\em The Riemann problem for Baer-Nunziato model of two phase flows}, J. Comput. Phys., 195: 434-464, 2004.

\bibitem{Bernetti} R. Bernetti, V.A. Titarev, E.F. Toro, {\em Exact solution of the Riemann problem for the shallow water equations with discontinuous bottom geometry}, J.Comput. Phys., 227: 3212-3243, 2008.

\bibitem{ChangHsiao3} T. Chang, L. Hsiao, {\em The Riemann problem and interaction of waves in gas dynamics}, Pitman Monographs, Longman Scientific and technica, 1989.

\bibitem{DalLeflochMurat4} G. Dal Maso, P.G. LeFloch, and F. Murat, {\em Definition and weak stability of nonconservative
products}, J. Math. Pures Appl., 74:  483-548, 1995.

\bibitem{Gallouet} T. Gallouet, J.M. Herard, N. Seguin. {\em Some approximation Godunov schemes to compute shallow-water equations with topography}, Comput. Fluids, 32: 889-899, 2003.

\bibitem{Goatin5} P. Goatin, P.G. LeFloch, {\em The Riemann problem for a class of resonant hyperbolic
systems of balance laws}, Ann. Inst. H. Poincare Anal. Non Lineaire, 21: 881-902, 2004.

\bibitem{Greenberg} J.M. Greenberg, A.Y. Leroux, {\em A well-balance scheme for the numerical processing of source terms in hyperbolic equations}, SIAM J. Numer. Anal., 33: 1-16, 1996.

\bibitem{Isaacson} E. Issacson, B. Temple, {\em Nonlinear resonance in systems of conservation laws}. SIAM J. Appl. Math., 52: 1260-1278, 1992.

\bibitem{Jin1} S. Jin, X. Wen, {\em An efficient method for computing hyperbolic systems with geometrical source terms having concentrations}, J. Comput. Math., 22:  230-249, 2004.

\bibitem{Jin2} S. Jin, X. Wen, {\em Two interface type numerical methods for computing hyperbolic systems with geometrical source terms having concentrations}, SIAM J. Sci. Comput., 26: 2079-2101, 2005.

\bibitem{LeflochThanh12} P.G. LeFloch, M.D. Thanh, {\em The Riemann problem for fluid flows in a nozzle with
discontinuous cross-section}, Commun. Math. Sci., 1: 763-797, 2003.

\bibitem{LeflochThanh13} P.G. LeFloch, M.D. Thanh, {\em The Riemann problem for shallow water equations with discontinuous topography}, Commun. Math. Sci., 5:  865-885, 2007.

\bibitem{LeflochThanh14} P.G. LeFloch, M.D. Thanh, {\em A Godunov-type method for shallow water equations with discontinuous topography in the resomant regime}, J. Comput. Phys., 230: 7631-7660, 2011.

\bibitem{LiTT} T.T. Li, W.C. Yu, {\em Boundary value problems of quasilinear hyperbolic system}, Duke University, USA, 1985.

\bibitem{Liu2} T.P. Liu, {\em Transonic gas flow in a duct of varying area}, Arch. Rational Mech. Anal., 80: 1-18, 1982.

\bibitem{SaurelAbgrall15} R. Saurel, R. Abgrall, {\em A multiphase Godunov method for compressible multifluid and
multiphase flows}, J. Comput. Phys., 150: 425-467, 1999.

\bibitem{Saurel17} R. Saurel, O. Lemetayer, {\em A multiphase model for compressible flows with interfaces, shocks, detonation waves and cavitation}, J. Fluid Mech., 431: 239-271, 2001.

\bibitem{Sheng} W.C. Sheng, Q.L. Zhang, {\em Interaction of the elementary waves of isentropic flow in a variable cross-section duct}. Commun. Math. Sci., 16:  1659-1684, 2018.

\bibitem{Smoller} J. Smoller, {\em Shock Waves and Reaction-Diffusion Equations}. Springer-Verlag, 1994.

\bibitem{Stoker} J.J. Stoker, {\em Water Waves}, Interscience, New York, 1957.

\bibitem{Thanh17} M.D. Thanh, {\em The Riemann problem for a nonisentropic fluid in a nozzle with discontinuous cross-sectional area}, SIAM J. Appl. Math., 69: 1501-1519, 2009.

\bibitem{Toro} E.F. Toro, {\em Shock capturing methods for free surface shallow water flows}, Wiley and Sons Ltd., 2001.

\bibitem{Wang} R.H. Wang, Z.Q. Wu, {\em Existence and uniqueness of solutions for some mixed initial boundary value problems of quasilinear hyperbolic systems in two independent variables} (in Chinese), Acta Scientiarum Naturalium of Jilin University 2: 459-502, 1963.

\bibitem{Zhang} Q.L. Zhang, {\em The interaction of elementary waves for nonisentropic flow in a variable cross-section duct}. To appear in Commun. Math. Sci..

\end{thebibliography}
 \end{document}